\documentclass[11pt]{article}
\usepackage[utf8]{inputenc}
\usepackage[a4paper, margin=0.7in]{geometry}
\usepackage[T1]{fontenc}
\usepackage{cite}
\usepackage{amsmath,amssymb,amsfonts}
\usepackage{algorithmic}
\usepackage{graphicx}
\usepackage{textcomp}
\usepackage{colortbl}
\definecolor{headergreen}{RGB}{0, 153, 0}
\definecolor{rowblue}{RGB}{220, 230, 241}
\usepackage{pgf-pie}
\usepackage{xcolor}
\usepackage{hyperref}
\usepackage{tikz}
\usetikzlibrary{arrows.meta, positioning}
\usepackage{amsmath}
\usepackage{venndiagram}
\usepackage[dvipsnames,table]{xcolor}
\definecolor{highlight}{RGB}{133,223, 255} 
\usepackage{multirow}
\usepackage{booktabs}
\title{A Systematic Review of Recent Advancements in PINN Augmented Deep Learning and Mathematical Modeling for Efficient Portfolio Management}

\author{
Bahadur Yadav,
Sanjay Kumar Mohanty \\
Department of Mathematics, School of Advanced Sciences \\
Vellore Institute of Technology, Vellore, Tamil Nadu 632014, India \\
Email: sanjaykumar.mohanty@vit.ac.in
}

\begin{document}
\date{}
\maketitle
\begin{abstract}
In finance, portfolio management is a traditional yet difficult problem that has drawn attention from practitioners and researchers for many years. However, there are still difficult technological problems that need to be solved. In the world of finance, managing a portfolio has never been easy. Selecting portfolios in a volatile market is made easier with the help of portfolio management. The goal of this review study is to present the concept of physics-informed neural networks because they provide a novel approach to directly incorporating physics and finance principles into the neural network's learning process. By doing so, physics-informed neural networks ensure that their forecasts are in line with established financial regulations and processes in addition to offering precise forecasts. Furthermore, this article provides an overview of the current state of research in portfolio optimization with the support of mathematical models, deep learning models and physics-informed neural networks. In addition, the advantages and disadvantages of various deep learning and mathematical modelling are discussed. Researchers and business professionals alike should find the data useful for advancing the field of investment management and trying out new portfolio management strategies. For this purpose, in this review work, emphasis is given to these factors. Finally, a few challenging issues and potential future directions are discussed, encouraging readers to consider fresh ideas in this field of study.
\end{abstract}


%
%
\thispagestyle{empty}


\section{Introduction}
Physics Informed Neural Networks (PINNs) incorporate physics concepts like differential equations and conservation laws into the design of the neural network. Because of their ability to integrate and learn from data while adhering to fundamental physical rules, PINNs are very helpful in solving challenges in finance, scientific computing,  biomedical engineering, physics-based simulations, etc.  Our focus is mostly on using PINNs to manage portfolios in the finance industry. Stock market prediction has become an enormous research topic in the field of research, and few research works have been accomplished with the help of PINNs. Asset management is a technique where time series prediction and portfolio optimization are the primary focus. In the world of finance, portfolio management involves managing the portfolio to minimize risk and optimize profit.

Portfolio management is a strategy used to handle one's wealth effectively, with the term portfolio typically referring to the collection of assets an individual owns. The practice has evolved over time, with distinct approaches in ancient and modern eras.
In ancient times, wealth primarily consisted of assets like land, agricultural produce, and various commodities. People managed their portfolios by exchanging these assets based on prevailing circumstances. Land and agriculture were the key sources of wealth, and individuals adjusted their holdings accordingly. This contrasts with modern portfolio management, which emphasizes systematic approaches tailored to individual risk tolerances. In today's world, portfolios typically comprise stocks, bonds and various securities. Managing modern portfolios involves optimizing returns while minimizing risk through asset allocation strategies.

In portfolio management, mathematical approaches give rigorous frameworks for risk quantification, portfolio allocation optimization and well-informed decision-making in the face of uncertainty. These strategies make use of statistical methods and mathematical models to improve portfolio performance and efficiently manage risk. Qualitative approaches provide insights into variables that are difficult to quantify, quantitative approaches give objectivity and rigor. Combining the two methods is often necessary for successful portfolio management in order to obtain a thorough grasp of investment opportunities and hazards. The integration of quantitative models with qualitative assessment facilitates informed decision-making by investors, enabling them to align their investment goals and risk tolerance.

Using machine learning and deep learning, we have been able to tackle numerous challenges related to the stock market application. Given the volatile nature of the market, managing a portfolio is always difficult in the finance sector. Moreover, maximizing return while minimizing risk is a difficult challenge. Numerous studies have been conducted on this.

In portfolio management, many research papers have been published on return prediction using machine learning and deep learning \cite{bali2023option} and portfolio models with return forecasting and transaction costs \cite{yu2020portfolio}. Robust active portfolio management \cite{erdogan2008robust} and a dynamic approach for the evaluation of portfolio selection performance under symmetric risk \cite{biglova2014portfolio}. Many papers have been proposed based on deep reinforcement learning for
trading \cite{zhang2020deep}. 

Deep learning techniques have been incredibly popular in recent years, and several deep learning techniques have been used in finance to predict the share market. Particularly when it comes to the topic of stock market prediction. Predicting the stock market is difficult when working with large amounts of data in finance, but deep learning can provide previously unimaginable results in this situation. Artificial Intelligence (AI) and mathematical modeling can be used in two different ways to analyze stock market predictions. Deep learning is an essential topic that frequently comes up in conversations about artificial intelligence.

\begin{figure*}
    \centering
    \includegraphics[width=15cm]{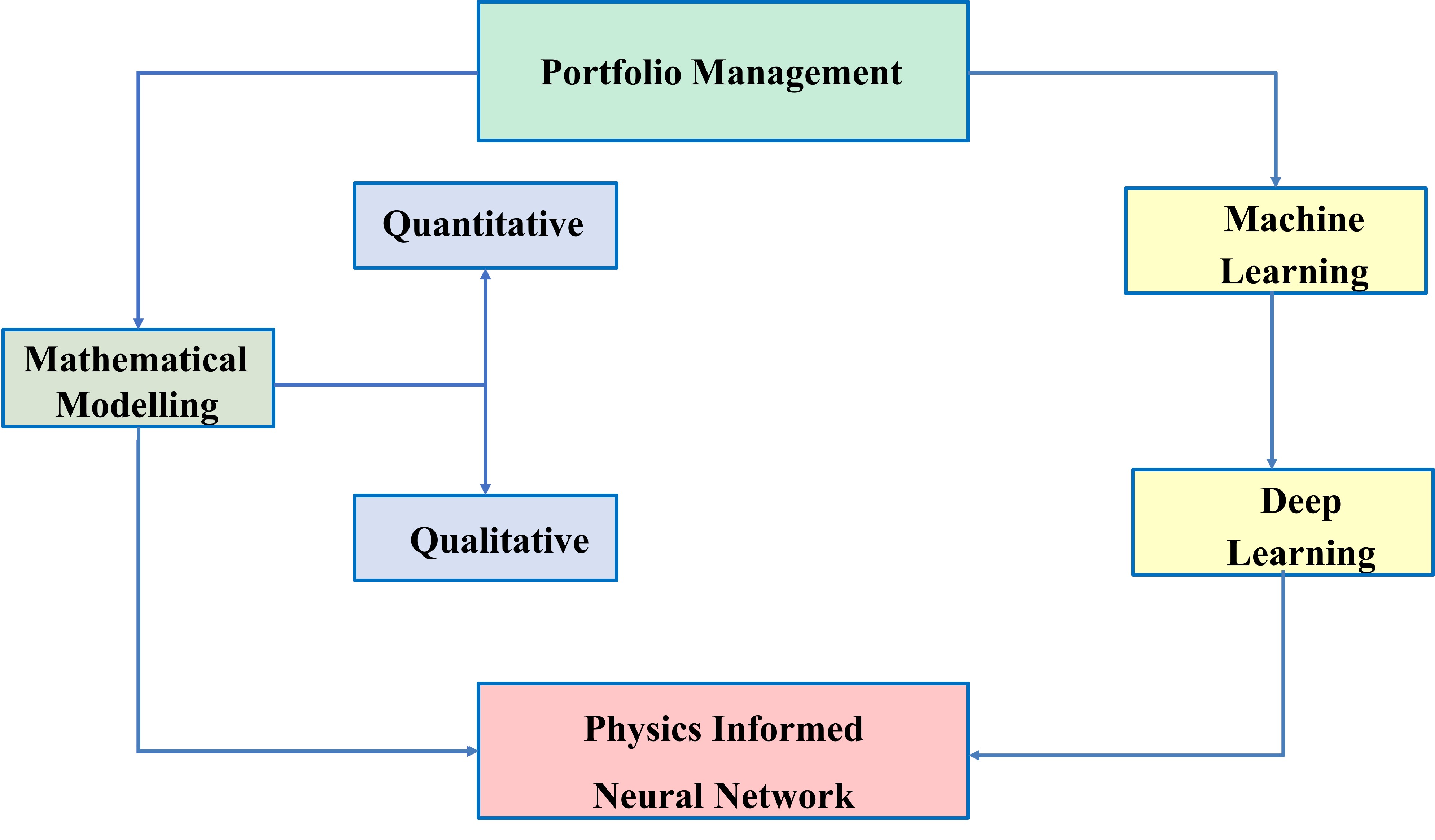}
    \caption{ Quick look at article organization generated by authors.}
    \label{fig:1}
\end{figure*}

Our primary tasks in this work are given as follows. Section \ref{1} provides the origin of portfolio management. Section \ref{2} focuses on mathematical methods for portfolio management. Section \ref{3} provides machine learning and deep learning for portfolio management. Next, physics-informed neural networks for portfolio management are discussed in section \ref{4}. In Section \ref{5}, a detailed discussion of performance measures and challenging issues is presented. Finally, a short conclusion ref{6} on the research gap and future directions is given. The
overview of this article is given in Figure \ref{fig:1}.

\section{Elements of Portfolio Management} \label{1}
A portfolio is an assortment of financial instruments, including futures and options, cash, bonds, stocks, and other assets, that are intended to yield returns in the future. Portfolio management is a strategy for managing one's wealth effectively. Both individuals and specialized financial organizations must manage their portfolios. For an individual or a company to grow, portfolio building and analysis are essential since they offer insights into possible investments and overall strategy.  The process of optimizing a portfolio is choosing the right combination of assets to meet predetermined goals, such as balancing risk or maximizing returns. Different optimization strategies are used by different individuals or companies to accomplish these goals effectively. In the upcoming subsection, we will provide a concise overview of the ancient and modern methodologies in portfolio management and examine their functionality during that era.

\subsection{Historical Developments}  
The practice has evolved, with distinct approaches of ancient and modern eras. The meaning of portfolio management in the past was completely different from what it is today. When it comes to diversification, land, commodities, livestock, and other assets made up the majority of their holdings at that time. These were the items that individuals traded according to their necessities.  Many challenges arose in the past when managing the portfolio. Due to a lack of knowledge, nothing was systemic. In terms of information availability, earlier portfolios were dependent on a small amount of data and low analytical standards to guide investment choices. They only did things based on necessity. Ancient Portfolio Theory(APT) failure might be explained by the fact that the absence of knowledge about technology created an environment in which nothing happened as planned, and things were exchanged based only on requirements. As a result of these shortcomings, modern portfolio theory emerged to supplant the outdated portfolio theories of the past.

Modern Portfolio Theory (MPT) primarily hinges on the examination of individual securities. It was produced by creating an ideal portfolio using an extremely basic and inadequate analytical method. It enhanced both investment practice and financial theory by focusing on the analysis of portfolio characteristics. 
Before the turn of the 20th century, APT entered the scene. It was suggested that it had a crucial role in the development of finance up until the 1952 introduction of Markowitz's mean-variance model. 
The primary reason APT was unable to provide a technical method for creating an ideal portfolio was that it lacked a mathematical procedure and statistical metrics to account for events that occurred in any statistically or mathematically meaningful way. It may be the primary cause of the ancient portfolio theory's failure. The absence of a basis in mathematics and statistics, as well as a variety of modeling and analysis techniques, all formed the foundation of modern portfolio theory, which was the field's earliest stage. Mathematical concepts then entered the scene, and MPT came into the picture and replaced the APT. The MPT provides the appropriate framework for market analysis, including the basis in mathematics and statistics, as well as several modeling techniques. We have now seen the transition of a portfolio from an ancient to a modern one. How both before and
after a mathematical model entered into the picture, it operated without a mathematical framework, making portfolio theory more effective than it was before. Having covered the disparities between ancient and modern portfolio theories, in the upcoming subsection, we will shift our focus to different types of portfolios, advanced techniques, and prominent journals in portfolio management.

\subsection{Exploring Diverse Portfolio Types: Foundations and Leading Journals in Portfolio Management}\label{sec22}
The portfolio is partitioned into five segments based on the level of risk and return. These segments are (i) Aggressive, (ii) Defensive, (iii) Income, (iv) Speculative, and (v) Hybrid. Each sort of portfolio has a unique set of attributes. The aggressive portfolio is primarily used to pursue higher returns. High returns are the goal of an aggressive portfolio, and in order to reach that goal, bigger risks are frequently taken. The enterprises in their early stages of growth are given priority. Conversely, the goal of a defensive portfolio is to reduce risk. Although it focuses on dividend earnings, an income portfolio is similar to a defensive portfolio in certain ways. The speculative portfolio is for high risk. The goal of the hybrid portfolio is to provide the best possible return at the lowest possible risk. Depending on the level of risk and return, it applies to various asset kinds.
Furthermore, based on the behaviour of the types of portfolios, it can be extended into four categories, viz., (i) Transaction cost portfolio, (ii) Robust portfolio, (iii) Regularized portfolio and (iv) Reinforcement learning portfolio.
Transaction fees or charges apply when we buy or sell assets during this time. The time penalization technique is used in the transaction cost portfolio to reduce transaction fees as much as possible. Conversely, the goal of a robust portfolio is to minimize the transaction costs during the sparse. An optimization process usually quantifies the estimates that are controlled within the regularized portfolio. Furthermore, a reinforcement learning portfolio handles the ongoing learning of market conditions necessary due to shifting market dynamics, mentioned in Tables \ref{tab:1} and \ref{tab:2}. Key elements, various approaches and best journals for portfolio management are listed in Tables \ref{tab:3}, \ref{tab:4}, and \ref{tab:5}.  Next, different mathematical approaches associated with portfolio management will be discussed for effective portfolio management.  
\begin{table*}[ht]
\centering
\caption{Types of portfolio}\label{tab:1}
\renewcommand{\arraystretch}{1.3}
\rowcolors{2}{blue!10}{blue!10} 
\begin{tabular}{|>{\bfseries}c|c|p{0.73\linewidth}|}
\rowcolor{green!30} 
\hline
\textbf{Sl. No.} & \textbf{Portfolio} & \textbf{Description}  \\
\hline
1 & Aggressive & Takes significant risks in an attempt to increase returns. \\
2 & Defensive & Minimizes risk and provides the least amount of return. \\
3 & Income & Quite similar to a defensive portfolio, but with an emphasis on dividend increases or other kinds of steady income. \\
4 & Speculative & Aims for a very high risk. Almost like gambling. \\
5 & Hybrid & Seeks to deliver the best possible return at the best possible level of risk. \\
\hline
\end{tabular}
\end{table*}


\begin{table*}[htbp]
\centering
\caption{Extended types of portfolio management}\label{tab:2}
\renewcommand{\arraystretch}{1.3}
\rowcolors{2}{blue!10}{blue!10} 
\begin{tabular}{|>{\columncolor{highlight}}c|c|p{0.64\linewidth}|}
\rowcolor{green!30} 
\hline
\textbf{Sl. No.} & \textbf{Portfolio} & \textbf{Description} \\
\hline
1 & Transaction Cost &  Seeks to lower transaction costs through the application of time penalty strategies. \\
\hline
2 & Robust &  Seeks to decrease transaction costs in a sparse and robust portfolio selection process. \\
\hline
3 & Regularized &  Seeks to minimize estimation error. \\
\hline
4 & Reinforcement Learning & Seeks to consistently learn about the shifting state of the market by utilizing a variety of strategies. \\
\hline
\end{tabular}
\end{table*}

\begin{table*}[htbp]
\centering
\caption{Key elements of portfolio management}\label{tab:3}
\renewcommand{\arraystretch}{1.3}
\rowcolors{2}{blue!10}{blue!10} 
\begin{tabular}{|>{\columncolor{highlight}}c|c|p{0.71\linewidth}|}
\rowcolor{green!30} 
\hline
\textbf{Sl. No.} & \textbf{Key element} & \textbf{Description} \\
\hline
1 & Asset Allocation &Appropriate asset allocation while considering both short and long-term goals is essential to efficient portfolio management. Because no two assets move the same, this demands a thorough grasp of the asset \\
\hline
2 & Diversification & One crucial component of an asset is how risk and return are distributed within it. It is challenging to have precise knowledge about any subset of an asset because of shifting market dynamics. As a result, diversification becomes a crucial component of portfolio management. \\
\hline
3 & Rebalancing & Rebalancing involves reallocating the assigned assets to the intended asset mix. Rebalancing is the practice of selling expensive assets and using the proceeds to purchase less expensive ones. \\
\hline
\end{tabular}
\end{table*}

\begin{table*}[htbp]
\centering
\caption{Exploring various approaches to portfolio management}\label{tab:4}
\renewcommand{\arraystretch}{1.3}
\rowcolors{2}{blue!10}{blue!10} 
\begin{tabular}{|>{\columncolor{highlight}}c|c|p{0.66\linewidth}|}
\rowcolor{green!30} 
\hline
\textbf{Sl. No.} & \textbf{Approach} & \textbf{Description} \\
\hline
1 & Passive Management & Long-term asset allocation is involved in this technique.
 \\
\hline
2 & Active Management & In order to obtain the highest return, assets must be actively purchased and sold on a regular basis. \\
\hline
\end{tabular}
\end{table*}

\begin{table*}
  \centering
\caption{List of the key  journals source for portfolio management}\label{tab:5}
\renewcommand{\arraystretch}{1.3}
\rowcolors{2}{blue!10}{blue!10} 
  \begin{tabular}{|>{\columncolor{highlight}}c|c|p{0.47\linewidth}|}
 \rowcolor{green!30} 
    \hline
\textbf{Sl. No.}  & \textbf {Journal name}  &  \textbf {Subject domain} \\
    \hline
 1 &   SIAM Journal On Control And Optimization   & Applied Mathematics \\
  \hline
 2 &    Journal of Finance   & Finance \\
  \hline
 3 &   The Review of Financial Studies  & Economics and Econometrics \\
  \hline
4 &    Review of Finance  & Financial Economics \\
 \hline
 5 &   Expert Systems with Applications  & Applied Mathematics \\
  \hline
6 &    Annals of Operation Research & Decision Sciences  \\

    \hline
  \end{tabular}
\end{table*}
\section{Mathematical Methods for Portfolio Management}\label{2}
 Qualitative and quantitative analyses are fundamental methods crucial for effective portfolio management, each serving distinct purposes. Qualitative analysis is crucial for validating the presence and distinctiveness, i.e., existence and uniqueness of solutions within mathematical models applied in portfolio management.  On the other hand, quantitative analysis is necessary to comprehend how various parameters influence the outcomes of these models. Table \ref{tab:6} explores advanced mathematical techniques essential for effective portfolio analysis and understanding of portfolio dynamics.
Next, we will discuss the qualitative and quantitative approaches to portfolio management and various techniques and strategies used to manage portfolios effectively.

\begin{table*}
\centering
\caption{ Application of  advanced mathematical techniques in finance(portfolio management).}\label{tab:6}
\renewcommand{\arraystretch}{1.3}
\rowcolors{2}{blue!10}{blue!10} 
\label{tab:stock_prediction}
\centering
\begin{tabular}{|>{\columncolor{highlight}}c|p{1.67cm}|p{6.1cm}|p{6.3cm}|}
\rowcolor{green!30} 
\hline
\textbf{Sl. No.}  & \textbf{Reference} & \textbf{Techniques} & \textbf{Advantage} \\
\hline
& &  &  \\
  1 & {\cite{hambly2017stochastic}} & Stochastic 
Volatility. &  Resilience to extreme events, enhanced stability and improved risk management. \\
\hline
& &  &  \\
 2 &  {\cite{ahmadi2019portfolio}} & Information \; Geometry and Entropy. & Capability to harness information flow and enhanced adaptability, as well as enhanced alpha generation. \\
\hline
 &  &  &\\
 3 &  {\cite{best1991sensitivity}} & Stochastic Calculus and Optimal \;Control. & Real-time adaptability, improved precision, and enhanced risk management. \\
\hline
 &  &  &\\
 4 &  {\cite{liang2013large}} & Large-Scale Optimization and Distributed Computing. & Scalability, speed, and efficiency, handling big data. \\
\hline
& & &\\
 5 & {\cite{davis1990portfolio}} & Optimization under Transaction Costs. & Realism in practical implementation, cost-efficient portfolio construction, and risk mitigation. \\
\hline
& & &\\
 6 &  {\cite{capponi2022systemic}} &  Risk Management and Dependence Modeling. & Enhanced diversification strategies, enhanced accuracy, comprehensive risk assessment, and improved tail risk management. \\
\hline
& & &\\
 7 & {\cite{zhang2020deeplear}} & Deep Learning and Alternative Data. & Optimize the portfolio through enhanced pattern recognition, improved prediction accuracy, and adaptability to changing market conditions. \\
\hline
& & &\\
 8 & {\cite{ban2018machine}} & Machine Learning and Artificial Intelligence. & Optimize the portfolio through enhanced insight generation,  unstructured data analysis, holistic decision-making, predictive capabilities, and data analysis power.\\
\hline
\end{tabular}
\end{table*}

\subsection{Qualitative Approach for Portfolio Management}
In a qualitative approach to portfolio management, concepts like stability, control sensitivity, existence, and uniqueness are essential for understanding the dynamics of portfolios and their management techniques. These things make it possible to find the best solutions in line with investors' goals and to develop flexible risk management and performance improvement tactics in the face of shifting market conditions.
 With an emphasis on European call option selling prices, Marcozzi \cite{marcozzi2008stochastic} examined stochastic optimum control for ultra-diffusion systems to improve portfolio management in a continuous time environment. On the other hand, Ekeland $\&$ Taflin \cite{ekeland2007optimal} created optimal bond portfolios. The goal of the author's work is to create a comprehensive framework for managing stocks and bonds in real time. Through the use of a Stochastic Partial Differential Equation (SPDE), Nadtochiy $\&$ Zariphopoulou \cite{nadtochiy2019optimal} provided a corresponding stochastic utility and explicit computing in the Black-Scholes model.
Hambly $\&$ Kolliopoulo {\cite{hambly2017stochastic}}  investigated a market model that includes assets that are susceptible to default, demonstrating that the empirical measure process has a huge portfolio limit. They addressed the intertemporal investment issues faced by mutual fund managers by suggesting a move away from high terminal fees and toward reduced interim fees.  In order to determine how sensitive a complex bond portfolio is to changes in the yield surface, Duedahl \cite{duedahl2016implementation} presented a method for computing stochastic duration. In order to beat a benchmark while adhering closely to its risk profile, Pesenti $\&$ Jaimungal \cite{pesenti2023portfolio} explored active portfolio management via lowering distortion risk measures. Expanding upon previous research, it presents different approaches inside a Wasserstein ball by utilizing isotonic projections in a broad market framework. Udeani $\&$  Ševčovič \cite{udeani2021application} examined the completely nonlinear Hamilton-Jacobi-Bellman (HJB)  equation for portfolio optimization. Using Banach's fixed-point theorem, solutions are achieved by converting it into a quasilinear parabolic equation, similar to porous media. Using dynamic programming and intertemporal manager rewards, Bensoussan et al. \cite{bensoussan2022inter} presented a novel approach to mutual fund management. It delivers original solutions, demonstrates adaptability in managing varying risk attitudes and recommends charge modifications to reduce client expenses without compromising management satisfaction. By applying variational inequalities to lump-sum and cash flow payoffs, Bensoussan et al. \cite{bensoussan2010real} proposed the best investment strategies in monopolistic and Stackelberg scenarios. Differential game interpretations emphasize the influence of competition while addressing issues that include fragmented markets and non-differentiable impediments. It combines utility maximization with choices about investments, portfolios, and consumption. Liu  {\cite{liu2007portfolio}} developed a portfolio selection in stochastic environments where the dynamic portfolio choice problem is specifically addressed and resolved by the author. In particular, when asset returns are quadratic, the method requires solving a system of ordinary differential equations. Acemoglu et al. {\cite{acemoglu2015systemic}} used systemic risk and stability in financial networks and presented a tractable theoretical framework for the study of the economic forces shaping the relationship between the structure of the financial network and systemic risk. Next, we will discuss the quantitative mathematical approach for portfolio management.
\subsection{Quantitative Approach for Portfolio Management }
In a quantitative approach to portfolio management, techniques like mathematical modelling,  statistical analysis, numerical methods, and optimization are essential for understanding and analysing financial data, managing risk, and constructing optimal portfolios. These factors help to improve investment decision-making processes and provide better results for investors.

The introduction of the first quantitative mathematical model developed by Markowitz {\cite{markowitz1952portfolio}}, which is known as the Mean-Variance (MV) model. He developed a quadratic programming model for picking a diversified stock portfolio, and in order to solve the problem, sophisticated nonlinear numerical techniques are used. Without ever resolving a mathematical programming issue, Elton et al. { \cite{elton1976simple}} introduced simple criteria for optimal portfolio selection and created decision criteria that enable one to arrive at the best solution to practical portfolio problems. Later, Konno $\&$ Suzuki {\cite{konno1995mean}} introduced a mean-variance skewness portfolio optimization model, which extends the classical Mean-Variance (MV) model. The skewness of the rate of returns of assets and the third-order derivatives of a utility function are significant factors in determining an optimal solution.  On the other hand, Konno $\&$ Yamazaki  {\cite{konno1991mean}} provided the Mean Absolute Deviation (MAD) from the mean as the risk measure. The model mean absolute deviation risk, which can overcome the challenges associated with Morkowitz’s classical model while retaining its advantages over equilibrium models. King  {\cite{king1993asymmetric}} formulated symmetric risk measures and tracking models for portfolio optimization under uncertainty and serves as an extension of Markowitz's MV model. An investigation by Best $\&$ Grauer {\cite{best1991sensitivity}} delivered sensitivity analysis for mean-variance portfolio problems. The primary objective of this study is to investigate the use of sensitivity analysis in mean-variance portfolio optimization through the use of a general parametric quadratic programming method.

Many authors proposed models for optimal portfolio management with transition costs. Pliska $\&$
Selby {\cite{pliska1994free}} initially introduced a model for optimal portfolio management with fixed transition costs, utilizing a geometric Brownian motion to describe asset prices. Later, Alexakis et al. {\cite{alexakis2007dynamic}} investigated a dynamic approach for the evaluation of portfolio performance under risk conditions. Lobo et al. {\cite{lobo2007portfolio}} presented portfolio optimization with linear and fixed transaction costs and addressed the difficulty of choosing a portfolio while taking transaction costs and risk exposure limits into consideration. Yu et al. {\cite{yu2020portfolio}} presented portfolio models with return forecasting and transaction costs. They investigated the efficacy of incorporating return forecasting into diverse portfolio models.

In order to optimize discounted consumption utility, Ahmadi Javid $\&$ Fallah Tafti{\cite{ahmadi2019portfolio}} developed Portfolio optimization with entropic value-at-risk. This study explores sample-based portfolio optimization utilizing the Entropic Value-at-Risk (EVaR) and introduces a primal-dual interior-point algorithm specifically designed for portfolio optimization incorporating EVaR. Erdogan et al.  {\cite{erdogan2008robust}} examined the robust active portfolio management
and provided solid models for managing an active portfolio in a market where there are transaction costs, controlling the effects of estimating mistakes in market parameter values on portfolio strategy performance is the aim of the robust models and allowing the imposition of additional side constraints like limits on the holdings in the portfolio. 

An additional investigation was conducted by  Kim et al. {\cite{kim2016time}} on the performance of time series momentum and volatility scaling. The Time Series Momentum (TSMOM) technique is reexamined by the authors in this research employing liquid futures contracts. Instead of a time series
momentum, they discover that volatility-scaling returns account for a major portion of their findings.
With consideration for fuzzy uncertainty,  Leon et al. {\cite{leon2000fuzzy}} presented a model for the problem of fuzzy mathematical programming for portfolio management. The purpose of this study is to employ fuzzy optimization strategies to effectively handle portfolios within the context of balancing risk and return. However, Biglova et al. {\cite{biglova2014portfolio}} introduced a portfolio selection problem in the presence of systemic risk and provided reward-risk measures that account for systemic risk and delivered the profitability of several strategies based on the forecasted evolution of returns and compared the optimal sample paths of future wealth obtained by performing reward-risk portfolio optimization on simulated data. Kral et al. {\cite{kral2019quantitative}} introduced a quantitative approach to project portfolio management and developed a mathematical framework for integer programming that takes bivalent variables into account in order to maximize project portfolios. Goetzmann $\&$ Kumar  {\cite{goetzmann2008equity}} developed equity portfolio diversification, where the study provides the portfolios of over 40,000 stock investment accounts from a sizable discount brokerage over six years in the recent history of the US capital market. Lin et al. {\cite{lin2023portfolio1}} discussed portfolio selection under systemic risk. This study provides a modified Sharpe ratio that can be used to create an ideal portfolio in the event of systemic problems. On the other hand, Alvarez et al. {\cite{alvarez2023optimal}} proposed the Almgren-Chriss model for brokerage contracts, optimizing trading for several clients trading a primary asset.  It can determine customer reservation values endogenously, allows the broker to choose clients strategically, and is computationally efficient when describing ideal portfolios. Miller $\&$ Yang {\cite{miller2017optimal}} achieved the ability to effectively handle time inconsistency in CVaR-based stochastic control and enable dynamic programming. Versatility is improved by extending to a variety of risk indicators, as demonstrated by portfolio optimization within CVaR limitations. Xidonas et al. {\cite{xidonas2018multiobjective}} proposed bridging mathematical theory with asset management practice. They studied a complete framework for portfolio optimization that closes the knowledge gap between complex mathematical theories. A comprehensive understanding of the various state-of-the-art methods is given in Table \ref{tab:7}. Next, we will discuss deep machine and deep learning methods for portfolio management and various techniques that can enhance portfolio management strategies.
\begin{table*}[htbp]
\centering
\caption{Portfolio optimization analysis models}\label{tab:7}
\renewcommand{\arraystretch}{1.3}
\rowcolors{2}{blue!10}{blue!10} 
\label{tab:portfolio_optimization}
\small 
\setlength{\tabcolsep}{4.4pt} 
\renewcommand{\arraystretch}{1.2} 
\begin{tabular}{|>{\columncolor{highlight}}c|p{2cm}|p{6.2cm}|p{6.3cm}|}
\rowcolor{green!30} 
\hline
\textbf{Sl. No.} & \textbf{References} & \textbf{Advantages} & \textbf{Limitations} \\
\hline
1 &  {\cite{markowitz2000mean}} & Used in situations where there are frequent decisions made about the portfolio and the return distribution is compact. & cannot be utilized in situations when the risk parameters are high. \\
\hline
2 &  {\cite{samuelson1970fundamental}} & When portfolio returns are not normally distributed, this technique can be used. & Optimization is limited to ensuring the best local solution. \\
\hline
3 &   {\cite{jorion1997value}} &  Relevant to a multitude of assets. &  Challenging to use in extensive portfolios. \\
\hline
4 &  {\cite{rockafellar2000optimization}} & Capable of employing a weighted average to handle severe loss. & There is no maximum loss that can be suffered according to CVaR. \\
\hline
5 & {\cite{konno1993mean}} &  Utilized when multivariate, regularly distributed asset returns are present. & May result in computation hazards from the covariance matrix being skipped. \\
\hline
6 &  {\cite{young1998minimax}} & Provides advantages that make sense when returns are not distributed properly. & Attentive to anomalies. Unusable in the absence of historical data. \\
\hline
7 &  {\cite{nawrocki1992characteristics,  brogan2008non}} & Fits with how investors view risk, punishes only negative deviations. & It takes time to calculate portfolio risk and alert about anomalies in data. \\
\hline
\end{tabular}
\end{table*}

\section{Machine Learning  and Deep Learning Approach in Portfolio Management}\label{3}
 In order to maximize portfolio performance, make better judgments, and adjust to shifting market conditions, machine learning provides several tools to portfolio managers. However, in practical investment scenarios, it is critical to properly validate models, control risks, and take into account the limitations of machine learning algorithms.  Deep learning approaches present effective methods to improve risk management, prediction, optimization, and decision-making, among other elements of portfolio management. Applying deep learning to actual investment settings requires addressing issues, including overfitting and low-quality data. Machine learning depends on feature engineering and simpler techniques, while deep learning uses multi-layered, complicated neural networks to enable autonomous hierarchical feature extraction. Next, we will discuss machine learning techniques for portfolio management.

\subsection{Recent Development of Machine Learning Technique for Portfolio
Management}
Machine learning is a vital tool for portfolio managers, helping them to maximize portfolio performance, make better judgments, and efficiently adjust to shifting market conditions. However, to fully utilize machine learning in portfolio management, it is imperative to recognize and overcome implementation-related obstacles, such as model interpretability and regulatory compliance. Bali et al.  {\cite{bali2023option}} proposed an option return predictability with machine learning and big data, and anticipated option returns using machine learning. It is observed that nonlinear models outperform and demonstrate that the complexity of machine learning models influences prediction. Christensen et al. {\cite{christensen2023machine}} anticipated a machine learning approach to volatility forecasting. The out-of-sample performance of several machine learning techniques for volatility forecasting is thoroughly analyzed and compared to the Heterogeneous Auto-Regressive (HAR) model.  Ma et al. {\cite{ma2021portfolio}} examined Portfolio optimization with return prediction using machine learning in order to outperform conventional time series models in portfolio creation. Papenbrock et al.  {\cite{author2021title}}  introduced synthetic correlations and explainable machine learning for constructing robust investment portfolios and used evolutionary algorithms to produce realistic correlation matrices. Using real-world data, Aithal et al.  {\cite{aithal2023real}} examined a real-time
portfolio management system utilizing machine learning techniques and discussed the key aspects of the entire process involved in portfolio selection, optimization, and management.  On the other hand, when applied to real data, Ban et al. {\cite{ban2018machine}} elaborated on a machine learning and portfolio optimization.  The authors investigated Performance Cross Validation (PBV) and Performance Regularization (PBR) in the portfolio optimization problem. To robustify the portfolio optimization, Lin $\&$ Taamouti {\cite{lin2023portfolio}} examined portfolio selection under non-Gaussianity and systemic risk. This research presents a novel performance ratio that accounts for both systemic risk and non-Gaussianity while constructing optimal portfolios. Gu et al. {\cite{gu2020empirical}} provided machine learning-based empirical asset pricing. By dynamically modifying asset weights in response to market conditions and risk objectives, Carta et al. {\cite{carta2021multi}} proposed that reinforcement learning approaches are utilized to optimize portfolio allocation strategies and possess the ability to weigh the trade-offs between risk and return, determining the best investment allocations to optimize portfolio returns and effectively manage risk. Deng et al. \cite{deng2016deep} determined that reinforcement learning has proven to outperform more conventional methods, such as supervised learning, in some financial settings. Complex patterns and adaptable tactics that may not be visible using traditional methods can be found by reinforcement learning algorithms, improving the results of portfolio management activities.  Gueant $\&$ Manziuk {\cite{gueant2019deep}} proposed deep reinforcement learning for market making in corporate bonds. Next, we will discuss deep learning techniques for portfolio management.
\subsection{Recent Development of Deep Learning Technique for Portfolio Management}
\begin{figure*}
    \centering
    \includegraphics[width=1\textwidth]{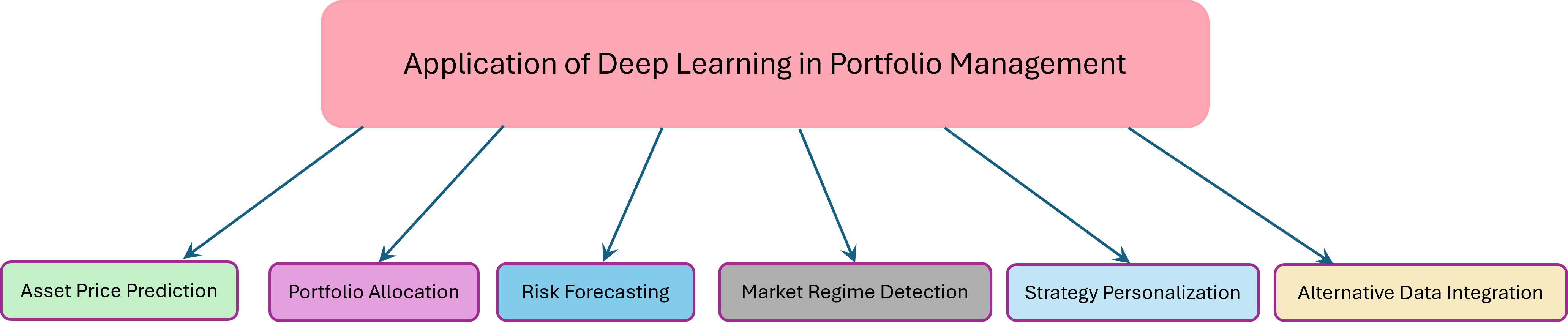}
    \caption{Deep learning application in portfolio management.}
    \label{10}
\end{figure*}

Deep learning has become increasingly popular in the modern era. Forecasting the stock market is an extremely challenging task. When applied to portfolio management, deep learning transforms decision-making by offering enhanced risk management, more precise forecasts, and flexible investment strategies that adjust to changing market conditions. Wang et al. {\cite{wang2020portfolio}} introduced portfolio formation with pre-selection using deep learning from long-term financial data. They combined Mean Variance (MV) models with Long Short-Term Memory (LSTM) networks to present a method for optimal portfolio building, with a focus on capturing long-term relationships in financial time series data. Soleymani $\&$
Paquet {\cite{soleymani2020financial}} discussed financial portfolio optimization with online deep reinforcement learning and restricted stacked autoencoder deep breath. Utilizing deep reinforcement learning and a constrained stacking autoencoder for feature selection. Yun et al. {\cite{yun2020portfolio}} introduced portfolio management via two-stage deep learning with a joint cost. They suggested using a joint cost function and a two-stage deep learning framework dubbed the Grouped-ETFs Model (GEM) to train a deep learning model for portfolio management. Zhu {\cite{zhu2020stock}} gave a Stock price prediction using the Recurrent Neural Networks (RNN) model. Here, a deep learning method for stock market forecasting that takes advantage of RNN benefits while working with time series data is discussed. Jiang et al.  {\cite{jiang2017deep}} classified a deep reinforcement learning framework for the financial portfolio management problem. They introduced a financial model, a free reinforcement learning framework designed to address portfolio management challenges. Zhang et al.  {\cite{zhang2020deeplear}} proposed deep learning for portfolio optimization. They used deep learning models to optimize the Sharpe ratio of a portfolio directly. Using gradient descent to update model parameters, this pipeline optimizes portfolio weights without requiring the typical forecasting step. Liang et al. {\cite{liang2018adversarial}} introduced an adversarial deep reinforcement learning in portfolio management. They incorporated three cutting-edge continuous reinforcement learning techniques into the realm of portfolio management. Also, they systematically evaluated the performance of these algorithms across diverse configurations, encompassing varying learning rates, objective functions, and feature combinations. Bertoluzzo $\&$ Corazza {\cite{bertoluzzo2012testing}} proposed testing different reinforcement learning configurations for financial trading. They developed and applied some automatic Financial Trading Systems (FTSs) based on differently configured reinforcement learning algorithms that optimize their behaviour in relation to the responses they get from the environment in which they operate without the need for a supervisor. Vo et al. {\cite{vo2019deep}} introduced a deep learning algorithm for decision-making and the optimization of socially responsible investments and portfolios. They presented a Deep Responsible Investment Portfolio (DRIP) model that uses a multivariate bidirectional LSTM neural network. Obeidat et al. {\cite{obeidat2018adaptive}} examined an adaptive portfolio asset allocation optimization with deep learning. They unveiled a neural network-based recommendation system for adaptive asset allocation. Zhang et al. {\cite{zhang2021universal}} provided an end-to-end framework that directly optimizes a portfolio by utilizing deep learning models and compared it to the classical two-step method. Chavan et al. {\cite{chavan2021intelligent}} introduced intelligent investment portfolio management using time-series analytics and deep reinforcement learning. This study suggests that Markov Decision Process (MDP) and reinforcement learning exhibit promise in developing practical models for stock trading
decisions. Minh et al. {\cite{minh2018deep}} proposed that deep learning techniques for stock pricing prediction enable more accurate forecasts and better risk management strategies in the face of market volatility. Fang et al. {\cite{fang2019research}} examined deep learning approaches to investing, such as high-frequency trading, option pricing, risk assessment, and transaction cost prediction, that help investors make better decisions and implement more effective portfolio management plans. Wang $\&$ Xu  {\cite{wang2018leveraging}} proposed an insurance business that can be of substantial benefit to portfolio management by implementing deep learning techniques, especially in the areas of fraud detection and risk assessment. Strengthening security protocols and lowering financial losses, Han et al. {\cite{han2018nextgen}} provided deep learning-based language technologies to augment anti-money laundering investigation for fraud detection in banking and online markets. Smalter Hall $\&$ Cook  {\cite{smalter2017macroeconomic}} employed macroeconomic indicator forecasting with deep neural networks to improve portfolio management by allowing for more dynamic portfolio adjustments, anticipating company swings, optimizing asset allocation strategies, facilitating more effective risk management, and reducing reliance on particular forecasting models. Loureiro et al. {\cite{loureiro2018exploring}} discussed deep neural networks for sales forecasting in fashion retail that improve portfolio management in retail markets through applications such as AR and sales forecasting models.
Fombellida et al. {\cite{fombellida2020tackling}} introduced a tackling business intelligence with
bioinspired deep learning that provides better skills for data analysis, more flexible learning, effective anomaly detection, sophisticated feature selection, and better decision-making processes to portfolio managers. Addo et al. {\cite{addo2018credit}} proposed  credit risk analysis using machine and deep learning
models that can improve portfolio management through improved loan pricing processes, more precise credit risk assessment, improved portfolio diversification, proactive risk mitigation strategies, efficient model performance evaluation, and improved decision-making processes. Go $\&$ Hong {\cite{go2019prediction}} proposed that pattern recognition and feature extraction from large, complex datasets are areas in which deep learning models shine. These skills allow portfolio managers to recognize asset correlations, market trends, and other significant variables affecting investment choices. Wang et al. {\cite{wang2020portfolio1}} presented a range of financial characteristics that can be forecasted by deep learning models, such as asset prices, volatility, and risk factors. These models improve forecast accuracy and help with risk management and portfolio optimization by utilizing past data and intricate relationships to learn. Song et al.{\cite{song2019study}} provided improved risk assessment and scenario analysis for portfolio management by capturing temporal interdependence and identifying trends in market behaviour. By determining the best asset allocations in accordance with risk-return objectives and restrictions, Tsang $\&$ Wong {\cite{tsang2020deep}} proposed an optimized deep learning algorithm for portfolio design. By integrating predictive models into portfolio optimization algorithms, s, Aggarwal $\&$ Aggarwal  {\cite{aggarwal2017deep}}  introduced deep learning techniques for financial market investment that make use of hierarchical decision models. Hu et al. {\cite{hu2021survey}} concluded that deep learning outperforms traditional approaches and allows for learning through interactions, which makes it a powerful tool and technique for portfolio management. It also optimizes trade execution and allocation tactics. Soleymani $\&$ Paquet{\cite{soleymani2021deep}} proposed deep learning models that can learn from fresh data and adjust to shifting market conditions. We have discussed different techniques which will help to manage the portfolio management, given in Tables \ref{tab:8} and \ref{tab:9}. Next, we will discuss physics-informed neural networks, which use mathematical models and deep learning models. The application of deep learning in portfolio management is presented 
in Figure \ref{10}.

\begin{table*}[htbp]
\centering
\caption{Comparative analysis of various state-of-the-art technique}\label{tab:8}
\renewcommand{\arraystretch}{1.3}
\rowcolors{2}{blue!10}{blue!10} 
\setlength{\tabcolsep}{4.4pt} 

\label{tab:state_of_art}
\small 
\setlength{\tabcolsep}{4pt} 
\renewcommand{\arraystretch}{1.1} 
\begin{tabular}{|>{\columncolor{highlight}}c|p{2.5cm}|p{5.6cm}|p{6.4cm}|}
\rowcolor{green!30} 
\hline
\textbf{Sl. No.} & \textbf{Methods} & \textbf{Advantages} & \textbf{Limitations} \\
\hline
1 & {\cite{markowitz1952portfolio}} & Determining  the trade off between risk and return. &  It is a reliance on assumptions.\\
\hline
2 & {\cite{konno1995mean}} &  Improves the model's capacity to identify the best solutions by taking higher order derivatives and skewness into account. & Might make things more difficult and need more computing power.\\
\hline
3 & {\cite{capponi2022systemic}} & Enhances risk management by customizing portfolio returns according to systemic occurrences. &  Complexity in figuring out the effects of systemic occurrences.\\
\hline
4 & {\cite{marcozzi2008stochastic}} & Provides a platform for managing stocks and bonds in real-time. &  The difficulty of putting stochastic control techniques into practice.\\
\hline
5 &   {\cite{jorion1997value}} & Straightforward and simple to put into practise. It functions well in cases where historical data is available, and the return is regularly distributed. & Useful in real-world situations, but challenging because returns are not always regularly distributed. \\
\hline
6 & {\cite{chaouki2020deep}} &  When past information is lacking, it can be used. & Losses might arise from initial allocation before profits are realized. \\
\hline
7 & {\cite{moody2001learning}} & Utilized in the process of optimizing performance metrics. & Not applicable in a situation requiring instant estimation.\\
\hline
8 &  {\cite{aboussalah2020continuous}} & In a multi-modal state space, it is necessary to use continuous action, which requires no forecasting of time series. &  Highly complex, the market needs to be in a situation today to make better predictions.\\
\hline
9 & {\cite{gatta2023meshless}} &  Effectiveness evaluation in the pricing of options. & There is a computational complexity.\\
\hline
10 & {\cite{bai2022application}} &  Improved convergence, stability, and accuracy compared to typical PINNs approaches. &  Strict validation is required across a range of datasets and contexts.\\
\hline
\end{tabular}
\end{table*}

\begin{table*}[htbp]
\centering
\caption{A Comparative analysis of traditional and creative techniques}\label{tab:9}
\renewcommand{\arraystretch}{1.3}
\rowcolors{2}{blue!10}{blue!10} 
\label{tab:comparative_study}
\small 
\setlength{\tabcolsep}{4pt} 
\renewcommand{\arraystretch}{1.1} 

\begin{tabular}{|>{\columncolor{highlight}}c|p{2.3cm}|p{5.7cm}|p{6.48cm}|}
\rowcolor{green!30} 
\hline
\textbf{Sl. No.} & \textbf{Techniques} & \textbf{Advantages} & \textbf{Drawbacks/Limitations} \\
\hline
1 & Statistical-based techniques & Simple to put into practice, suitable for shorter time spans and smaller data volumes, and simple to understand. & Disregards the significance of several features. \\
\hline
2 & Regression-Based techniques & Easy to understand, simple to apply, construable, and able to use past data to learn repeatedly. & Efficiency is greatly influenced by the selected learning rate.\\
\hline
3 & Bayesian-based techniques & Easily expandable, in contrast to other machine learning-based methods, suitable for challenges involving reasoning. & There is no one-size-fits-all way to build networks out of data. Not every training pattern will be employed during the training process. \\
\hline
4 & Neural Network Based techniques & Superior learning model in comparison to other conventional methods. A distributed memory architecture, fault tolerance, and the ability to function with limited knowledge. & When working with big datasets, this should be utilized. Lack of explicit guidelines and inability to be self-explanatory in determining the network structure. \\
\hline
5 & Reinforcement-Based techniques & Able to resolve challenging optimization issues, the learning patterns of people and learning models are extremely similar and have the potential to learn adaptively in response to shifting conditions. & It should not be applied to tiny data sets and straightforward problems. Over-rewarding can have a negative impact on learning. The learning process may be constrained by the dimensionality curse. \\
\hline
6 & Evolutionary-Based techniques &With its straightforward concept and application of non-linear constraints and non-stationary situations, it surpasses traditional performance methods. Strong resistance to dynamic fluctuations and the capacity for parallelism. & More computing power is required. \\
\hline
7 & Quantum-based techniques & Qubits help solve complicated and large-scale computing problems since they can exponentially expand storage. Quicker than any alternative approach. & Comparing quantum computers to conventional computers, the energy requirement is substantially higher. \\
\hline
\end{tabular}
\end{table*}

\subsection{Explainable Artificial Intelligence in Portfolio Management}
\section{Physics Informed Neural Network in Portfolio Management}\label{4}
A specific type of neural network designed to integrate physical concepts or equations into the learning phase is known as Physics-Informed Neural Networks (PINNs). PINNs combine domain-specific knowledge and deep learning techniques to tackle different physics problems effectively. PINNs are specifically designed to tackle supervised learning assignments aimed at simulating physical processes and systems. This includes a wide range of applications, such as the modeling of quantum mechanics, fluid flows, and other fields where complex physical principles regulate the underlying processes. 

\begin{figure}[h!]
    \centering
    \begin{tikzpicture}
        \begin{scope}[blend mode=overlay]
            \fill[blue!30] (-2,0) circle (3);  
            \fill[red!30] (2,0) circle (3);   
        \end{scope}

        \draw[green!60!black] (-2,0) circle (3);
        \draw[green!60!black] (2,0) circle (3);

        \node at (-3,0) {\textbf{Neural Network}};
        \node at (2.5,0) {\textbf{Physics}};
        \node at (0,0) {\shortstack{\textbf{Physics} \\ \textbf{Informed} \\ \textbf{Neural} \\ \textbf{Network}}};
    \end{tikzpicture}
    \caption{Simple architecture of physics-informed neural networks generated by authors}
    \label{Fig:2}
\end{figure}
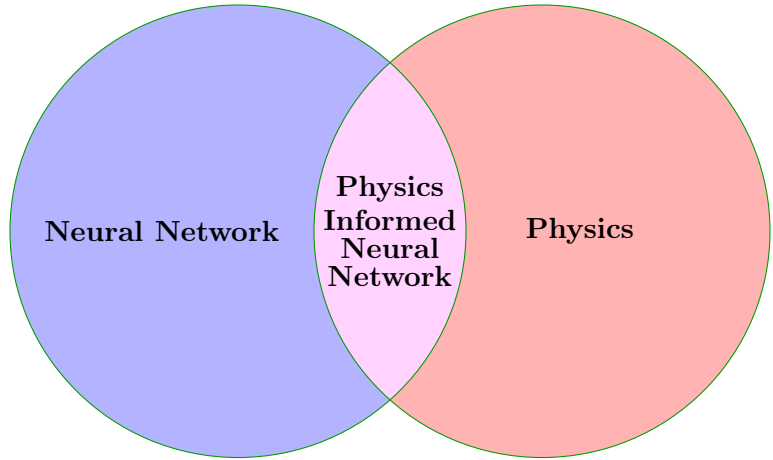

In Figure {\ref{Fig:2}}, PINNs are hybrids between well-established physics concepts and data-driven supervised neural networks. PINNs provide a special benefit by combining data-driven learning with strict adherence to established physics equations. This allows them to build models that efficiently utilize available data while also adhering to underlying physical rules. Because of this method, PINNs can generate more robust models with less data needed and still allow for precise extrapolation outside of the dataset that is now available. In a financial environment, this is becoming more complex and data-rich, and PINNs provide an effective tool for portfolio managers to make decisions. PINNs can assist in ways that traditional approaches might miss when it comes to risk management, investment opportunity identification, and portfolio performance optimization by utilizing AI and machine learning techniques. PINNs might improve decision-making.

\subsection{ Introduction to Physics Informed Neural Networks }
Differential equations and other pertinent physical concepts are integrated into the training process of a customized neural network defined as a PINN. By employing well-established physical rules to direct the neural network, this method enhances the model's precision and capacity for generalization—particularly when working with sparse training data.
PINNs basically use a neural network to approximate a function that satisfies a particular differential equation. In the course of training, the network aims to decrease the error between its forecasts and the given data, as well as the error between its forecasts and the preset physical principles represented by the differential equations. A loss function and a neural network architecture are combined to create PINNs. A physics-informed regularization term and a data mismatch term make up the two main parts of the loss function. These components work together to provide the foundation for teaching the PINNs to successfully strike a balance between precision and adherence to fundamental physical laws.

A revolutionary method for directly incorporating physics and financial regulations into the neural network learning process is presented by physics-informed neural networks. Not only do PINNs produce accurate forecasts, but PINNs also guarantee that the projections comply with relevant financial rules and guidelines. Neural network (NN) approaches are revolutionized by PINNs, which integrate model equations such as Partial Differential Equations (PDEs) into their architecture. This novel method has shown to be quite helpful in solving a variety of problems, including PDEs, fractional equations, integral differential equations, and stochastic PDEs. In addition to forcing neural networks (NNs) to fit observed data, PINNs serve as an adaptable multi-task learning framework that simultaneously minimizes PDE residuals. These PINNs, which are specially made to address PDE-based challenges, show promise as effective instruments in scientific machine learning. PINNs are essentially based on neural network training techniques that minimize a complete loss function to approximate PDE solutions. This function includes other important elements like the beginning and boundary conditions over the space-time domain and the PDE residuals found at particular collocation points in the domain. Fundamental to PINNs, deep-learning networks are trained to anticipate solutions at certain points in the integration domain of a differential equation. PINNs are unique in that they incorporate a residual network that contains the physics equations that regulate the system, which is a new feature. PINNs use an unsupervised training process, which does not require labelled data from experimental or simulation runs in the past. As a mesh-free method, the PINNs algorithm essentially transforms the direct solution of governing equations into an optimization problem with a loss function at its core. By adding a residual term that is taken from the governing equation, this method reinforces the loss function by integrating the mathematical model deeply into the neural network. Refinement and shrinking of the admissible solution space are achieved by this residual term acting as a restricting element.

PINNs are a class of machine learning models that combine traditional physics-based modeling with neural networks to solve forward and inverse problems in physics. Regression analysis is used to look for unknown parameters in a PDEs system by analyzing experimental data and PDE solutions. This process is known as the inverse problem. By simultaneously improving their set of neural network functions to supply unknown values at the selected set of collocation locations and fulfil the empirically measured values acting as the PDEs solutions, PINNs can tackle inverse issues.  Moreover, the losses incurred in the forward problem are validated along with the losses derived from the experimental data. To obtain the lowest possible losses, the settings of the DNNs and PDEs are tuned jointly. Different paths to tailored model extensions or applications are promised by the PINNs technique. Multiple unknown parameter determinations in inverse issues are easily conceivable. The capacity of PINNs to quickly and readily incorporate outside data, such as measurement or simulation data, into the total loss function is one of its primary characteristics.
Here is a brief explanation.

\begin{figure*}
    \centering
    \includegraphics[width=0.9\textwidth]{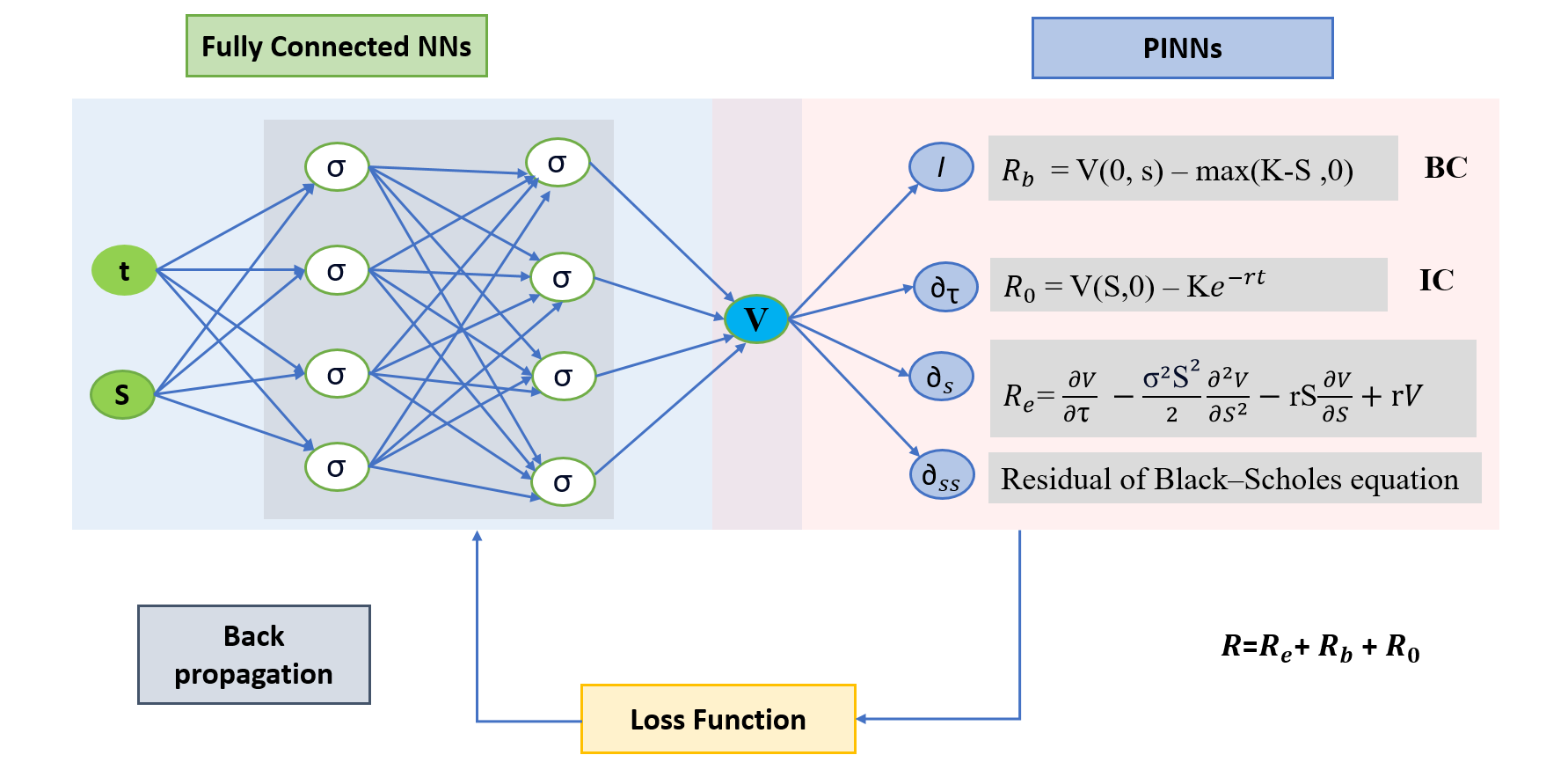}
    \caption{Overview of PINNs, schematic of PINNs framework. A fully connected neural network is used to approximate the solution $V(S, \tau)$, which is then applied to construct the residual
loss $R_e$, boundary conditions loss $R_b$ and initial conditions loss $R_0$.}
    \label{Fig:3}
\end{figure*}

\begin{enumerate}
\item \textbf{Forward Problem}: In a forward problem, we have a physical system described by differential equations. The goal is to predict the behaviour of the system over time or space, given certain initial and boundary conditions. Consider the heat conduction equation describing how temperature changes over time in a material:
\[ \frac{\partial k}{\partial t} = \alpha \nabla^2 k \]
where \( k \) is temperature, \( t \) is time, and \( \beta \) is the thermal diffusivity. Given the initial temperature distribution \( k(x,0) \) and boundary conditions, we want to predict \( k(x,t) \).
    
\item \textbf{Inverse Problem}: In an inverse problem, using the system's observations, we expect to determine the underlying physics model's parameters or initial/boundary conditions. Given temperature measurements at certain points in a material at various times, we want to infer the thermal diffusivity \( \beta \) and the initial temperature distribution \( k(x,0) \).
\end{enumerate}
By jointly solving the forward and inverse problems, PINNs can provide accurate solutions even when data is limited or noisy, making them valuable tools in many scientific and engineering applications.
In the heat conduction problem, we would train a neural network to approximate the temperature field \( k(x,t) \) while also adjusting the parameters such as \( \beta \) to match the observed data. This allows us to predict both the evolution of temperature over time and infer the properties of the material simultaneously.
We have discussed a comparative analysis between the inverse and forward methods within PINNs for portfolio management, mentioned in Table \ref{tab:10}. \\
\begin{table*}
    \centering
   \caption{Comparison of inverse and forward methods in PINN for portfolio management}\label{tab:10}
   \renewcommand{\arraystretch}{1.3}
\rowcolors{2}{blue!10}{blue!10} 
    \label{tab:pinn_comparison}
    \begin{tabular}{|>{\columncolor{highlight}}c|p{2.2cm}|p{3.8cm}|p{3.10cm}|p{4.2cm}|}
   \rowcolor{green!30} 
        \hline
        \textbf{Sl. No.} & \textbf{Aspect} & \textbf{Inverse Method  } & \textbf{Forward Method } & \textbf{ Application} \\
        \hline
        1 & \textbf{Description} & Determining parameters based on observed results. & Predicting results from known parameters. & Manage investment portfolios using PINN's ideas. \\
        \hline
        2 & \textbf{Functionality} & To estimate parameters, understand system dynamics. & Forecast future results scenarios. & Evaluate risk, optimize allocation, and estimate parameters. \\
        \hline
        3 & \textbf{Risk \;Modeling} & Analyze the risk-return of various assets and portfolios. & To estimate risk, simulate various market situations. & Recognize the risk involved with investing techniques. \\
        \hline
        4 & \textbf{Optimization} & Determine optimal portfolio allocations. & Combine assets in an optimal way to achieve desired goals. & Minimize risk and maximize returns. \\
        \hline
        5 & \textbf{Scenario Analysis} & Play with various market conditions. & Examine market fluctuations on portfolio performance. & Recognize the possible results under different circumstances. \\
        \hline
        6 & \textbf{Parameter Estimation} & Calculate parameters like expected returns, volatilities, etc. & Estimate parameters based on past data. &  Identify essential asset characteristics for decision making. \\
        \hline
    \end{tabular}
\end{table*}
In the context of PINNs, the solutions of standard Black–Scholes equation  PDEs may be approximated using a neural network, which is stated as
\begin{equation}
\frac{\partial V}{\partial t} + \frac{1}{2} \sigma^2 S^2 \frac{\partial^2 V}{\partial S^2} + rS\frac{\partial V}{\partial S} - rV   = 0,\text{ for } (S, t) \in (0,\infty) \times (0,T)
\end{equation}

\begin{equation}
V(S,0) = \max(K - S, 0);
\end{equation}
\begin{equation}
V(0, \tau) = Ke^{r\tau}; \quad V(L, \tau) = 0;
\end{equation}
{where  } K { is the strike price}; and  L { is some suitably large value by truncating the original } $(0,\infty)$ { domain to } (0, L]. $r$ represents the risk-free interest rate, and $\sigma > 0$ denotes the volatility; $V(S, \tau)$ represents the fair value of a vanilla European option at time $\tau$ with asset value $S$ at that time; $T$ denotes the expiration date of the contract, $\tau = T - t$ denotes the time to maturity, and $t \in [0, T]$ is the instantaneous time, mentioned in Figure \ref{Fig:3}.

\subsection{ Physics Informed Neural Network in Finance}
\begin{figure*}
    \centering
    \includegraphics[width=0.9\textwidth]{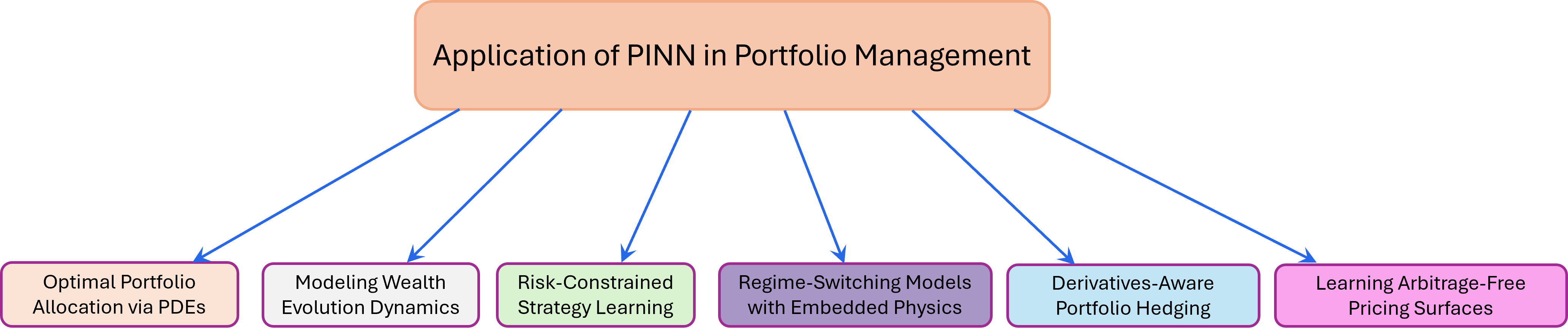}
    \caption{PINN application in portfolio management.}
    \label{Fig:6}
\end{figure*}

To address pricing issues for American multi-asset options, Gatta et al. { \cite{gatta2023meshless}} presented PINNs, a revolutionary deep learning technique. Using neural networks to solve intricate partial differential equations, it improves convergence for free boundary problems and changes conventional approaches. It creates a parametric model for additional market analysis and tests the efficacy of PINNs in option pricing. For the purpose of predicting volatility in financial markets, Kim et al. {\cite{kim2024physics}} proposed a unique architecture that combines convolutional transformers with neural networks inspired by physics. It provides better volatility surface estimation and performs better in numerical tests than current techniques. When analytical solutions to the Black-Scholes equation are hard to come by, Wang et al.{\cite{wang2023deep}} turned to PINNs to address option pricing. PINNs, which draw inspiration from Raissi et al. {\cite{raissi2019physics}}, provide a simplified method that is as accurate and computationally efficient as standard numerical approaches, indicating that they could be useful for the pricing of complicated options. In order to improve performance, Bai et al. { \cite{bai2022application}} demonstrated an improved version of the PINNs technique that incorporates local adaptive activation functions. It has been effectively implemented in finance models, including the Black-Scholes and Ivancevic option pricing models. When compared to typical PINN approaches, it exhibits better convergence, stability, and accuracy, and it advances differential equation solutions in a variety of domains. The application of PINN in portfolio management is presented 
in Figure \ref{Fig:6}.

\section{ DISCUSSION}\label{5}
\begin{figure}
    \begin{center}
\begin{tikzpicture}
\pie[
    text=legend,
    color={blue!50!white, orange, green!70!black, cyan!80},
    sum=auto,
    before number=\bfseries,
    after number=\%
]{
    35/Mathematical Modeling,
    30/Deep Learning,
    20/Portfolio Management,
    15/Conventional Machine Learning
}
\end{tikzpicture}
\end{center}
    \caption{ Distribution of research focus areas.}
    \label{fig:4}
\end{figure}
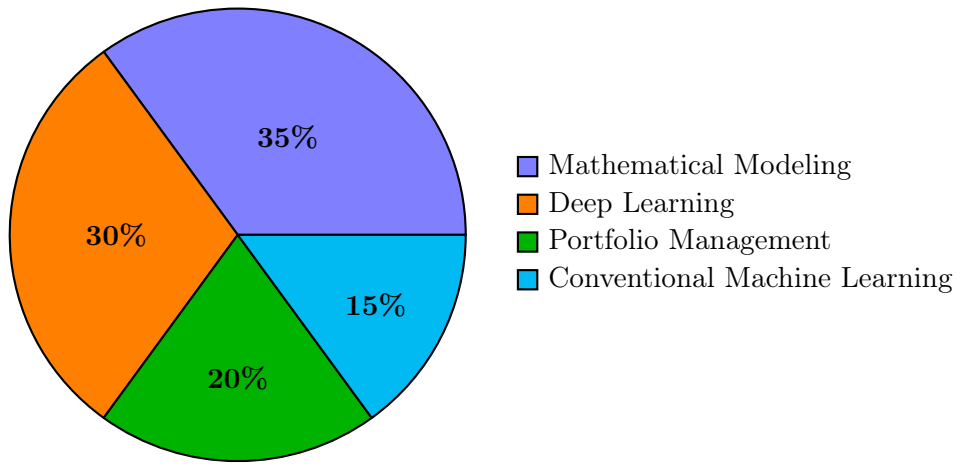

\begin{figure*}
    \centering
    \includegraphics[width=18cm]{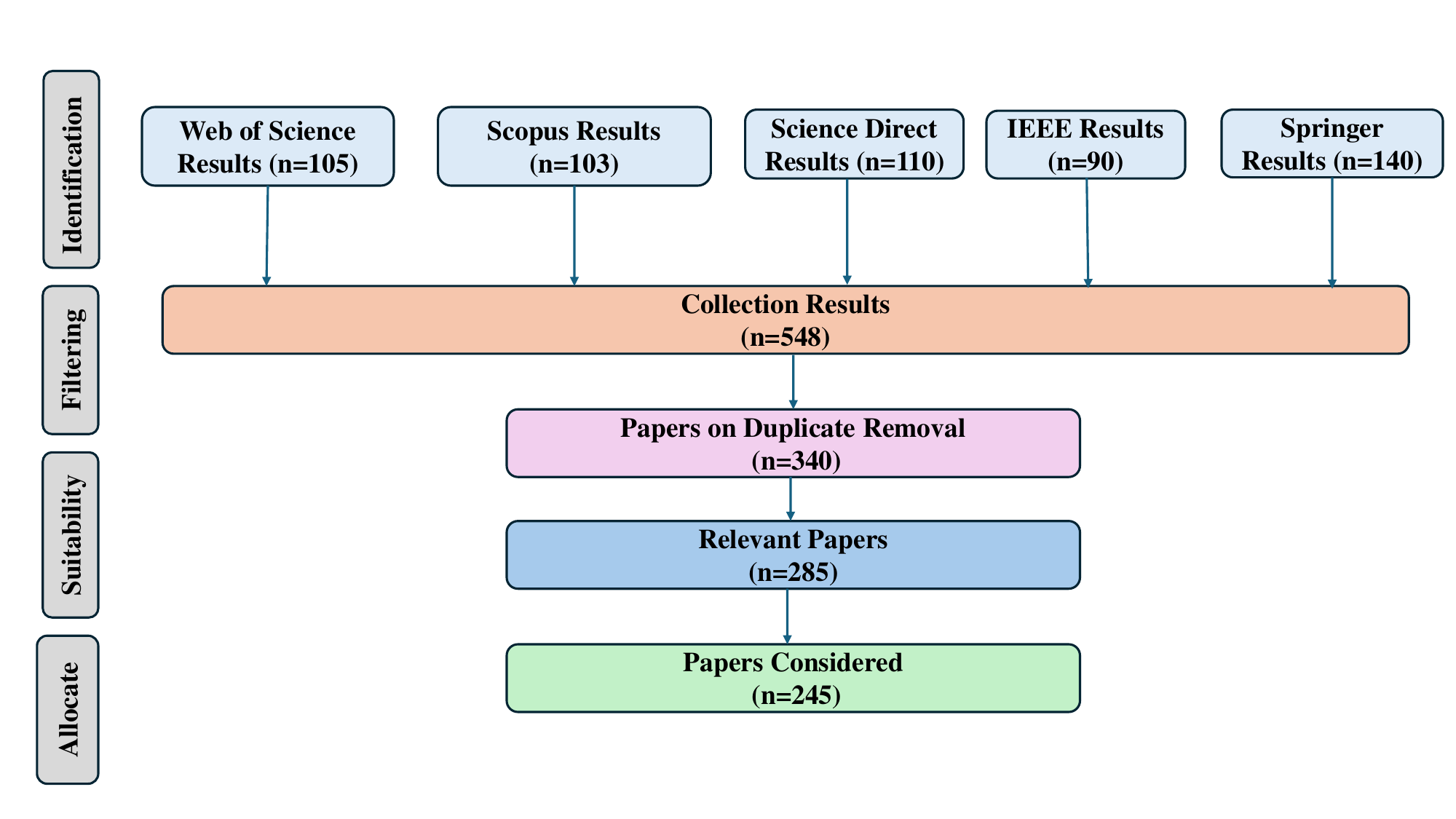}
    \caption{  Search strategy for portfolio management with a mathematical model and deep learning.}
    \label{fig:5}
\end{figure*}

This section covers an overview of performance measures, challenging issues, and future scopes, which are covered in detail in the following three subsections. Additionally, the distribution of research focus areas
and the search strategy of articles is highlighted in
Figures \ref{fig:4} and \ref{fig:5}, respectively.

\subsection{Performance Measures in Portfolio Management}
In portfolio management, classification evaluation metrics, regression evaluation metrics, and performance ratios are crucial tools that each has a unique role in evaluating and improving investment strategies. Classification models are used to identify and classify investment risks, while regression evaluation metrics are applied to continuous prediction tasks. Furthermore, by calculating risk-adjusted returns, comparing outcomes to benchmarks, and evaluating overall investment efficacy, performance ratios are used to analyze how effectively portfolio strategies perform. These tools work together to give investors and managers a thorough understanding of portfolio performance, enabling them to make informed decisions and keep risk and return at an ideal level. All of these topics are discussed in further detail in the next subsections that follow.

\subsubsection{Classification Metrics in Portfolio Management}
In finance, classification metrics are essential for making decisions about things such as investment risk classification, fraud detection, and the prediction of credit defaults.  For example, these metrics classify financial situations into distinct classifications, such as identifying whether a borrower is likely to default or whether a transaction is fraudulent. Accurate classification in fraud detection helps financial organizations detect suspicious activity in real-time, reducing losses and improving security. In credit scoring, metrics are used to determine whether a consumer is a high-risk or low-risk borrower, which affects interest rate allocations and loan approvals.  Likewise, investment risk models categorize assets or portfolios according to risk levels to facilitate capital allocation and portfolio development. Key measures obtained from the confusion matrix—accuracy, precision, recall, F1-score, and AUC-ROC—are used to assess how effective these classification models are.  These metrics offer information about a model's overall accuracy as well as its capacity to identify infrequent but crucial occurrences like fraud or defaults. We have discussed different evaluation metrics that help to manage portfolio performance, as mentioned in Table \ref{tab:11}.

\begin{table}[h!]
\centering
\caption{Evaluation metrics for classification models in portfolio management}\label{tab:11}
\renewcommand{\arraystretch}{1.3}

\rowcolors{2}{blue!10}{blue!10} 
\begin{tabular}{>{\bfseries}p{2.5cm} p{4cm} p{4.5cm} p{4.5cm}}
\rowcolor{green!30} 
\toprule
\textbf{Metric} & \textbf{Application} & \textbf{Advantages} & \textbf{Drawbacks/Limitations} \\
\midrule
Accuracy & Credit score and stock trend forecasting. & Straightforward and understandable measure. & For unbalanced datasets, it may be deceptive. \\
Precision & Fraud detection and alerts about insider trading. & Low false positive rate, prevents needless warnings. & Might overlook actual positives. \\
Recall & Predicting defaults and identifying bankruptcy. & Identifies more true positives. & May result in more false positives. \\
F1-Score & Classification of loan risk and churn forecast. & Useful for balancing recall and precision. & Not as interpretable as individual metrics. \\
AUC-ROC & Comparing models and risk rankings. & Robust against imbalance and independent of threshold. & Doesn't show the ideal threshold for decisions. \\
\bottomrule
\end{tabular}
\end{table}

\subsubsection{ Regression Metrics in Portfolio Management}
In financial modeling and time series forecasting, regression evaluation metrics such as {Mean Absolute Error (MAE)}, {Mean Squared Error (MSE)}, {Root Mean Squared Error (RMSE)}, {Mean Absolute Percentage Error (MAPE)}, {Symmetric MAPE (sMAPE)}, and the {Coefficient of Determination (R\textsuperscript{2})} play a crucial role in assessing the accuracy and reliability of predictive models. These metrics measure how much the true observed values (such as profits, stock prices, and NAVs) differ from the predictions made by the model.  While it is easier to understand and less susceptible to outliers, MAE gives a simple average of absolute errors.  MSE and RMSE are especially helpful in risk-sensitive financial applications like credit score and volatility predictions since they penalize greater errors more severely. By normalizing the error as a percentage, MAPE and sMAPE make it possible to compare data across scales, which is particularly useful when working with several assets or portfolios. R\textsuperscript{2} provides information on the percentage of variation that the model explains, which helps in performance comparison. Depending on the financial environment, each metric has a unique role, and analysts and portfolio managers must choose the best one based on whether they value interpretability, resilience to outliers, or comparable analysis across assets. We have discussed different evaluation metrics that help manage portfolio performance, as summarized in Table \ref{tab:12}.



\begin{table}[h!]
\centering
\caption{Evaluation metrics for regression models in portfolio management}\label{tab:12}
\renewcommand{\arraystretch}{1.3}

\rowcolors{2}{blue!10}{blue!10} 
\begin{tabular}{>{\bfseries}p{3cm} p{4.5cm} p{4.5cm} p{4.5cm}}
\rowcolor{green!30} 
\toprule
\textbf{Metric} & \textbf{Application} & \textbf{Advantages} & \textbf{Drawbacks/Limitations} \\
\midrule
MAE & Forecasting portfolio returns and NAV. & Interpretable and less impacted by outliers. & Disregards the direction of errors and treats them all equally. \\
MSE & Risk assessment and volatility forecasting. & Penalizes major mistakes, which is helpful in situations where risk is vital. & Outlier sensitivity, squared units. \\
RMSE & Forecasting asset prices and predicting returns. & Penalizes major errors in the same units as the data. & More difficult to comprehend than MAE and sensitive to outliers. \\
MAPE & Predicting income or returns (across assets) & Scale-independent and based on percentages & Unable to deal with zero values, skewed toward tiny denominators. \\
sMAPE & Making predictions about low-value assets or pricing. & Fair handling of excessive forecasts. & Less reasonable. \\
R\textsuperscript{2} & Goodness-of-fit for the model. & Variance explanation, commonly used benchmark. & May be misleading when using non-linear models. \\
\bottomrule
\end{tabular}
\end{table}
\subsubsection{Performance Ratios in Portfolio Management}
Performance ratios are crucial tools in portfolio management for assessing risk-adjusted returns and directing investment choices.  Investors can learn how successfully profits offset volatility by using the Sharpe Ratio, which calculates the excess return per unit of total risk. When evaluating portfolios with asymmetric return distributions, the Sortino Ratio is more appropriate since it avoids penalizing upside volatility by only taking into account downside risk.  The Treynor Ratio is appropriate for portfolios that are well-diversified because it employs beta as the risk metric and concentrates on systematic risk. For active managers looking to beat an index, the Information Ratio is especially pertinent since it compares portfolio returns to a benchmark after accounting for tracking inaccuracy.  In contrast to beta, which measures a portfolio's susceptibility to market fluctuations and helps identify exposure to systematic risk, alpha measures a portfolio's excess return in relation to its projected return based on beta, indicating management skill. Furthermore, by comparing return to the maximum drawdown, the Calmar Ratio provides information on risk-adjusted performance over longer timeframes.  By using these ratios, portfolio managers can evaluate risk exposure, compare methods, and choose the best portfolios to suit investors' risk-return preferences.  We have discussed different evaluation metrics that help manage portfolio performance, as summarized in Table \ref{tab:13}.

\begin{table}[h!]
\centering
\caption{Evaluation of performance ratios in portfolio management}\label{tab:13}
\renewcommand{\arraystretch}{1.3}

\rowcolors{2}{blue!10}{blue!10} 
\begin{tabular}{>{\bfseries}p{3cm} p{4.5cm} p{4.5cm} p{4.5cm}}
\rowcolor{green!30} 
\toprule
\textbf{Ratio} & \textbf{Application} & \textbf{Advantages} & \textbf{Drawbacks/Limitations} \\
\midrule
Sharpe Ratio & Examine portfolio performance in relation to overall risk. & Calculable and broadly accepted. & Anticipates normal returns and sanctions upward volatility equally. \\
Sortino Ratio & Consider downside risk solely while evaluating performance accounting. & Emphasizes damaging volatility and is more practical for returns that are not typical. & Downside deviation is less frequently utilized and demands more data. \\
Treynor Ratio & Evaluate performance in the context of systematic and market risk. & Beneficial for diversified portfolios. & Depends on the stability of beta. \\
Information Ratio & Assess the skill of the active manager in comparison to a standard. & Accounts for tracking errors. & Needs a suitable benchmark and is sensitive to the choice of benchmark. \\
Alpha & Calculates extra return over anticipated return. & Shows the manager's added value and is easy to understand. & Could be deceptive if market expectations are off. \\
Beta & Evaluates the sensitivity of a portfolio to market risk. & Aids in the building of portfolios and risk attribution. & Disregarding unsystematic risks. \\
Calmar Ratio & Returns are compared to the largest historical drawdown. & focuses on protecting capital and is appropriate for long-term performance. & Requires extensive historical data and is impacted by traumatic prior occurrences. \\
\bottomrule
\end{tabular}
\end{table}

\subsection{Challenging issues}
The following difficulties will arise in building a mathematical model for portfolio management and incorporating the deep learning approach:
\begin{enumerate}
    \item It is difficult to better comprehend asset behavior and portfolio management techniques when theoretical frameworks are integrated with continuously changing market dynamics and financial knowledge.

    \item It is difficult to create a mathematical model that takes into consideration a portfolio's macro-level framework, such as sectoral movements, market volatility, and risk exposure, in order to facilitate more precise investment planning.

    \item It is challenging to create a deep learning model for portfolio management from scratch since it requires an extensive amount of diverse high-quality financial data to provide efficient training and guarantee robust generalization capabilities.

    \item Model performance may suffer from problems like class imbalance in financial datasets (unusual market crashes vs. typical times) and multi-modal objectives (e.g., balancing return, risk, and liquidity).  Furthermore, diversified market conditions require adaptable strategies, which are challenging to develop and execute with little contextual knowledge.
\end{enumerate}

\subsection{Future Scopes}
This section provides an overview of current trends, gaps in the literature, and potential avenues for future research in order to provide fresh perspectives for additional study.
The review demonstrates that over the past decade, portfolio issues related to optimization have been successfully formulated using mathematical modelling, deep learning models, and PINNs. It might be feasible to employ various multi-objective optimization strategies in future studies, including compromise programming. Furthermore, the model was included in the study on multi-objective portfolio optimization, which took into account actual trading limitations. In light of these limitations, it would appear sensible and beneficial to highlight the multi-objective portfolio optimization problem. Another noteworthy finding is that, when addressing portfolio optimization problems with real-world restrictions, many of the examined research have opted to use heuristics as opposed to exact solution approaches. Effective mathematical modeling and deep learning models for portfolio management have been proposed by several researchers, as was previously indicated. Here are some potential future trends that could emerge

\begin{figure*}
    \centering
    \includegraphics[width=13cm]{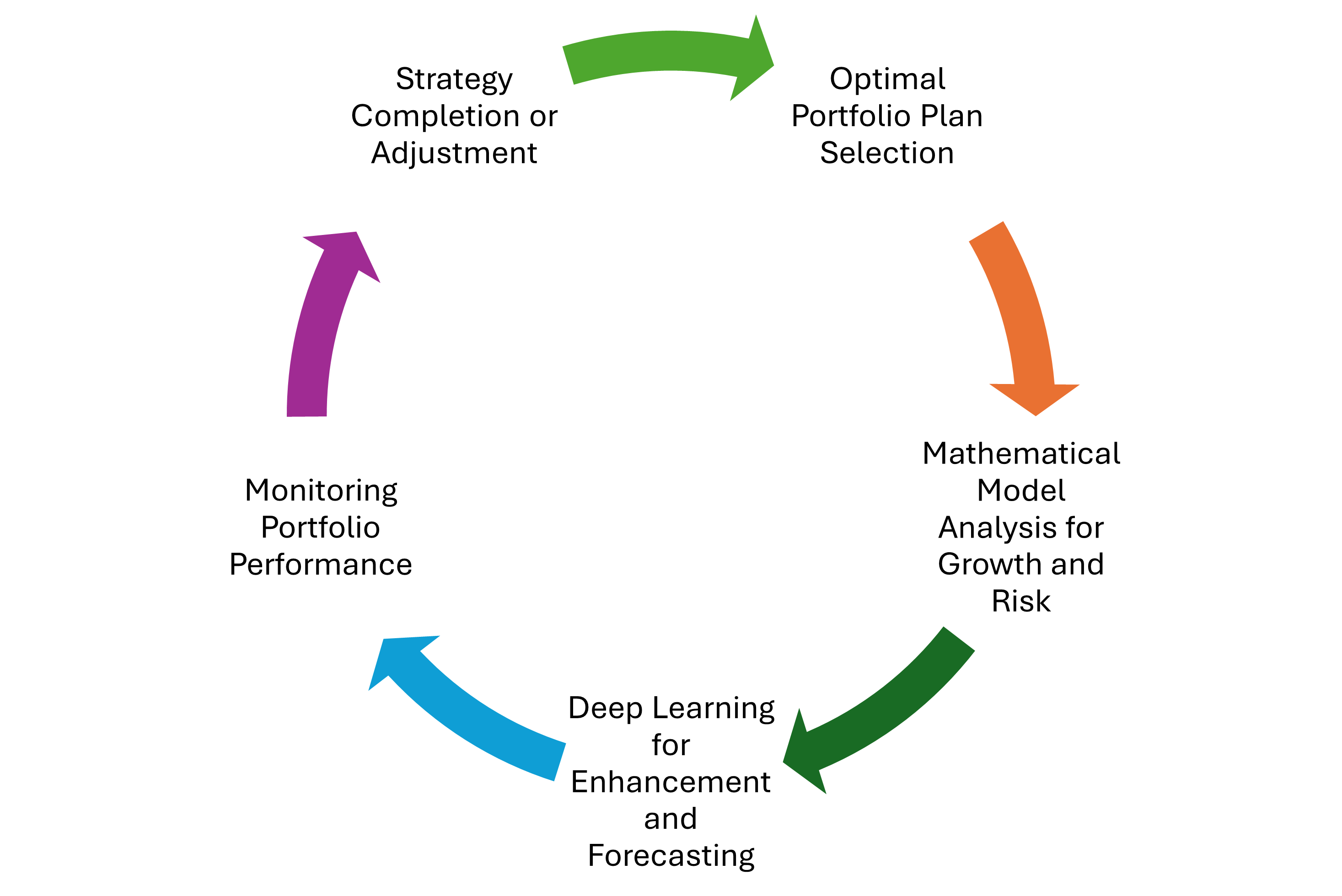}
    \caption{Portfolio management strategy with mathematical model and deep learning.}
    \label{fig:6}
\end{figure*}

\begin{itemize}
  \item To create deep learning, PINNs, and mathematical models that are therapeutically useful while taking into account the most recent advancements in finance.
\item To build models of portfolio management that are more effective by integrating the advantages of deep learning, PINNs, mathematical modeling, and expert knowledge.
\item To create integrated frameworks that use deep learning methods, PINNs, and mathematical models to comprehend portfolio management from various angles, perhaps helping investors in the stock market.

\item With the advantages of fuzzy optimization, a future trend in portfolio management will be to use hybrid models of fuzzy and PINNs techniques to solve complex multi-objective optimization problems.
\item Further research is necessary to find reliable and effective optimization methods that can produce better optimization scenarios for portfolio management.
\item To create integrated frameworks that combine deep learning techniques with mathematical models to analyze financial market dynamics and portfolio performance under various investment strategies. This could help portfolio managers make the best decisions and manage risk, as shown in Figure \ref{fig:6}.
\end{itemize}

\section{Conclusion}\label{6}
In this review article, we introduce the recently developed PINN technique and discuss its objectives, methodologies, and limitations. This article provided an extensive overview of the most current developments in deep learning, PINNs, and mathematical modeling for portfolio management. We have examined several deep learning, PINNs, and mathematical modeling approaches that are currently used in portfolio management. Further, we discussed several portfolio kinds that are utilized for managing and optimizing portfolios. This article presented a comparison of several traditional methods for portfolio optimization, and the application of mathematical and deep learning programming techniques to portfolio optimization problems has been extensively reviewed in the literature. Numerous portfolio optimization techniques were looked into, and works in this area were compared and evaluated. Additionally, this summary has highlighted a few possible areas for further investigation. Finally, future prospects and difficult challenges for this field of study have been discussed. We believe that this article will help researchers to analyze the recent developments in portfolio management.

\section*{CONFLICT OF INTEREST}
The authors of this article have disclosed that they have no potential conflicts of interest about the research, writing, or publication of the work.

\section*{Acknowledgements}

S. K. Mohanty acknowledges partial support from the Department of Science and Technology, Government of India, under grant number SR/FST/MS-II/2023/139-VIT Vellore.

\section*{Author contributions statement}

B.Y. and S.K.M. developed the review's structure and scope. After conducting the literature search, B.Y. arranged the studies that were gathered. Critical insights and a thorough examination of the findings were supplied by S.K.M. Both authors helped write and edit the work, and they gave their approval for submission.

\section*{Data availability statement}
This study has not generated any new data.

\bibliographystyle{IEEEtran}
\bibliography{sample}

@article{konno1995mean,
  title={A mean-variance-skewness portfolio optimization model},
  author={Konno, Hiroshi and Suzuki, Ken-ichi},
  journal={Journal of the Operations Research Society of Japan},
  volume={38},
  number={2},
  pages={173--187},
  year={1995},
  publisher={The Operations Research Society of Japan}
}

@article{alexakis2007dynamic,
  title={A Dynamic approach for the evaluation of portfolio performance under risk conditions},
  author={Alexakis, Christos and Balios, Dimitris and Stavraki, Sophia and others},
  journal={Invest. Manag. Financ. Innov},
  volume={4},
  pages={16--24},
  year={2007}
}

@article{king1993asymmetric,
  title={Asymmetric risk measures and tracking models for portfolio optimization under uncertainty},
  author={King, Alan J},
  journal={Annals of Operations Research},
  volume={45},
  pages={165--177},
  year={1993},
  publisher={Springer}
}

@article{pliska1994free,
  title={On a free boundary problem that arises in portfolio management},
  author={Pliska, Stanley R and Selby, Michael J},
  journal={Philosophical Transactions of the Royal Society of London. Series A: Physical and Engineering Sciences},
  volume={347},
  number={1684},
  pages={555--561},
  year={1994},
  publisher={The Royal Society London}
}

@article{xidonas2018multiobjective,
  title={Multiobjective portfolio optimization: bridging mathematical theory with asset management practice},
  author={Xidonas, Panos and Hassapis, Christis and Mavrotas, George and Staikouras, Christos and Zopounidis, Constantin},
  journal={Annals of Operations Research},
  volume={267},
  pages={585--606},
  year={2018},
  publisher={Springer}
}

@article{lobo2007portfolio,
  title={Portfolio optimization with linear and fixed transaction costs},
  author={Lobo, Miguel Sousa and Fazel, Maryam and Boyd, Stephen},
  journal={Annals of Operations Research},
  volume={152},
  pages={341--365},
  year={2007},
  publisher={Springer}
}

@article{liu2007portfolio,
  title={Portfolio selection in stochastic environments},
  author={Liu, Jun},
  journal={The Review of Financial Studies},
  volume={20},
  number={1},
  pages={1--39},
  year={2007},
  publisher={Oxford University Press}
}

@article{kral2019quantitative,
  title={Quantitative approach to project portfolio management: proposal for Slovak companies},
  author={Kral, Pavol and Valjaskova, Viera and Janoskova, Katarina},
  journal={Oeconomia Copernicana},
  volume={10},
  number={4},
  pages={797--814},
  year={2019}
}

@article{konno1991mean,
  title={Mean-absolute deviation portfolio optimization model and its applications to Tokyo stock market},
  author={Konno, Hiroshi and Yamazaki, Hiroaki},
  journal={Management science},
  volume={37},
  number={5},
  pages={519--531},
  year={1991},
  publisher={INFORMS}
}

@article{best1991sensitivity,
  title={Sensitivity analysis for mean-variance portfolio problems},
  author={Best, Michael J and Grauer, Robert R},
  journal={Management science},
  volume={37},
  number={8},
  pages={980--989},
  year={1991},
  publisher={INFORMS}
}

@article{wang2020portfolio,
  title={Portfolio formation with preselection using deep learning from long-term financial data},
  author={Wang, Wuyu and Li, Weizi and Zhang, Ning and Liu, Kecheng},
  journal={Expert Systems with Applications},
  volume={143},
  pages={113042},
  year={2020},
  publisher={Elsevier}
}

@article{ma2021portfolio,
  title={Portfolio optimization with return prediction using deep learning and machine learning},
  author={Ma, Yilin and Han, Ruizhu and Wang, Weizhong},
  journal={Expert Systems with Applications},
  volume={165},
  pages={113973},
  year={2021},
  publisher={Elsevier}
}

@misc{liang2018adversarial,
      title={Adversarial Deep Reinforcement Learning in Portfolio Management}, 
      author={Zhipeng Liang and Hao Chen and Junhao Zhu and Kangkang Jiang and Yanran Li},
      year={2018},
      eprint={1808.09940},
      archivePrefix={arXiv},
      primaryClass={q-fin.PM}
}

@article{vo2019deep,
  title={Deep learning for decision making and the optimization of socially responsible investments and portfolio},
  author={Vo, Nhi NY and He, Xuezhong and Liu, Shaowu and Xu, Guandong},
  journal={Decision Support Systems},
  volume={124},
  pages={113097},
  year={2019},
  publisher={Elsevier}
}

@article{soleymani2020financial,
  title={Financial portfolio optimization with online deep reinforcement learning and restricted stacked autoencoder—DeepBreath},
  author={Soleymani, Farzan and Paquet, Eric},
  journal={Expert Systems with Applications},
  volume={156},
  pages={113456},
  year={2020},
  publisher={Elsevier}
}

@article{yun2020portfolio,
  title={Portfolio management via two-stage deep learning with a joint cost},
  author={Yun, Hyungbin and Lee, Minhyeok and Kang, Yeong Seon and Seok, Junhee},
  journal={Expert Systems with Applications},
  volume={143},
  pages={113041},
  year={2020},
  publisher={Elsevier}
}

@article{aithal2023real,
  title={Real-Time Portfolio Management System Utilizing Machine Learning Techniques},
  author={Aithal, Prakash K and Geetha, M and Dinesh, U and Savitha, Basri and Menon, Parthiv},
  journal={IEEE Access},
  volume={11},
  pages={32595--32608},
  year={2023},
  publisher={IEEE}
}

@inproceedings{zhu2020stock,
  title={Stock price prediction using the RNN model},
  author={Zhu, Yongqiong},
  booktitle={Journal of Physics: Conference Series},
  volume={1650},
  number={3},
  pages={032103},
  year={2020},
  organization={IOP Publishing}
}

@article{obeidat2018adaptive,
  title={Adaptive portfolio asset allocation optimization with deep learning},
  author={Obeidat, Samer and Shapiro, Daniel and Lemay, Mathieu and MacPherson, Mary Kate and Bolic, Miodrag},
  journal={International Journal on Advances in Intelligent Systems},
  volume={11},
  number={1},
  pages={25--34},
  year={2018}
}

@article{ban2018machine,
  title={Machine learning and portfolio optimization},
  author={Ban, Gah-Yi and El Karoui, Noureddine and Lim, Andrew EB},
  journal={Management Science},
  volume={64},
  number={3},
  pages={1136--1154},
  year={2018},
  publisher={INFORMS}
}

@inproceedings{leon2000fuzzy,
  title={Fuzzy mathematical programming for portfolio management},
  author={Leon, Teresa and Liern, Vicente and Vercher, Enriqueta},
  booktitle={Financial Modelling},
  pages={241--256},
  year={2000},
  organization={Springer}
}

@article{yu2020portfolio,
  title={Portfolio models with return forecasting and transaction costs},
  author={Yu, Jing Rung and Chiou, Wan Jiun Paul and Lee, Wen Yi and Lin, Shun Ji},
  journal={International Review of Economics \& Finance},
  volume={66},
  pages={118--130},
  year={2020},
  publisher={Elsevier}
}

@article{ahmadi2019portfolio,
  title={Portfolio optimization with entropic value-at-risk},
  author={Ahmadi Javid, Amir and Fallah Tafti, Malihe},
  journal={European Journal of Operational Research},
  volume={279},
  number={1},
  pages={225--241},
  year={2019},
  publisher={Elsevier}
}

@article{zhang2020deeplear,
  title={Deep learning for portfolio optimization},
  author={Zhang, Zihao and Zohren, Stefan and Roberts, Stephen},
  journal={The Journal of Financial Data Science},
  volume={2},
  pages={8--20},
  year={2020},
  publisher={Institutional Investor Journals Umbrella}
}

@article{zhang2021universal,
  title={A universal end-to-end approach to portfolio optimization via deep learning},
  author={Zhang, Chao and Zhang, Zihao and Cucuringu, Mihai and Zohren, Stefan},
  journal={arXiv preprint arXiv:2111.09170},
  year={2021}
}

@article{chavan2021intelligent,
  title={Intelligent Investment Portfolio Management using Time-Series Analytics and Deep Reinforcement Learning},
  author={Chavan, Sachin and Kumar, Pradeep and Gianelle, Tom},
  journal={SMU Data Science Review},
  volume={5},
  number={2},
  pages={7},
  year={2021}
}

@article{jiang2017deep,
  title={A deep reinforcement learning framework for the financial portfolio management problem},
  author={Jiang, Zhengyao and Xu, Dixing and Liang, Jinjun},
  journal={arXiv preprint arXiv:1706.10059},
  year={2017}
}

@article{kim2016time,
  title={Time series momentum and volatility scaling},
  author={Kim, Abby Y and Tse, Yiuman and Wald, John K},
  journal={Journal of Financial Markets},
  volume={30},
  pages={103--124},
  year={2016},
  publisher={Elsevier}
}

@article{elton1976simple,
  title={Simple criteria for optimal portfolio selection},
  author={Elton, Edwin J and Gruber, Martin J and Padberg, Manfred W},
  journal={The Journal of finance},
  volume={31},
  number={5},
  pages={1341--1357},
  year={1976},
  publisher={JSTOR}}

@article{lin2023portfolio,
  title={Portfolio selection under non-gaussianity and systemic risk: A machine learning based forecasting approach},
  author={Lin, Weidong and Taamouti, Abderrahim},
  journal={International Journal of Forecasting},
  year={2023},
  publisher={Elsevier}
}

@article{zhang2020deep,
  title={Deep reinforcement learning for trading},
  author={Zhang, Zihao and Zohren, Stefan and Stephen, Roberts},
  journal={The Journal of Financial Data Science},
  volume={2},
  number={2},
  pages={25--40},
  year={2020}
}

@article{erdogan2008robust,
  title={Robust active portfolio management},
  author={Erdogan, Emre and Goldfarb, Donald G and Iyengar, Garud},
  journal={Journal of Computational Finance},
  volume={11},
  number={4},
  pages={71},
  year={2008}
}

@article{bertoluzzo2012testing,
  title={Testing different reinforcement learning configurations for financial trading: Introduction and applications},
  author={Bertoluzzo, Francesco and Corazza, Marco},
  journal={Procedia Economics and Finance},
  volume={3},
  pages={68--77},
  year={2012},
  publisher={Elsevier}
}

@article{capponi2022systemic,
  title={Systemic risk-driven portfolio selection},
  author={Capponi, Agostino and Rubtsov, Alexey},
  journal={Operations Research},
  volume={70},
  number={3},
  pages={1598--1612},
  year={2022},
  publisher={INFORMS}
}

@article{goetzmann2008equity,
  title={Equity portfolio diversification},
  author={Goetzmann, William N and Kumar, Alok},
  journal={Review of Finance},
  volume={12},
  number={3},
  pages={433--463},
  year={2008},
  publisher={Oxford University Press}
}

@article{acemoglu2015systemic,
  title={Systemic risk and stability in financial networks},
  author={Acemoglu, Daron and Ozdaglar, Asuman and Tahbaz-Salehi, Alireza},
  journal={American Economic Review},
  volume={105},
  number={2},
  pages={564--608},
  year={2015},
  publisher={American Economic Association 2014 Broadway, Suite 305, Nashville, TN 37203}
}

@article{biglova2014portfolio,
  title={Portfolio selection in the presence of systemic risk},
  author={Biglova, Almira and Ortobelli, Sergio and Fabozzi, Frank J},
  journal={Journal of Asset Management},
  volume={15},
  number={5},
  pages={285--299},
  year={2014},
  publisher={Springer}
}

@article{bali2023option,
  title={Option return predictability with machine learning and big data},
  author={Bali, Turan G and Beckmeyer, Heiner and Moerke, Mathis and Weigert, Florian},
  journal={The Review of Financial Studies},
  volume={36},
  number={9},
  pages={3548--3602},
  year={2023},
  publisher={Oxford University Press}
}

@article{author2021title,
 title= {Matrix Evolutions: Synthetic Correlations and Explainable Machine Learning for Constructing Robust Investment Portfolios},
  author= {Papenbrock, J. and Schwendner, P. and Jaeger, M. and Krügel, S.},
  journal= {The Journal of Financial Data Science},
  volume= {3},
  number= {2},
  pages= {51--69},
  year= {2021},
  publisher={Springer}
}

@article {lin2023portfolio1 ,
	title = {Portfolio Selection under Systemic Risk},
        author = {Lin, Weidong and Olmo, Jose and Taamouti, Abderrahim},
	eissn = {1538-4616},
	issn = {0022-2879},
	journal = {Journal of Money, Credit and Banking},
	publicationstatus = {In Press},
	publisher = {Wiley},
	year = {2024},
}

@article{christensen2023machine,
  title={A machine learning approach to volatility forecasting},
  author={Christensen, Kim and Siggaard, Mathias and Veliyev, Bezirgen},
  journal={Journal of Financial Econometrics},
  volume={21},
  number={5},
  pages={1680--1727},
  year={2023},
  publisher={Oxford University Press}
}

@article{gu2020empirical,
  title={Empirical asset pricing via machine learning},
  author={Gu, Shihao and Kelly, Bryan and Xiu, Dacheng},
  journal={The Review of Financial Studies},
  volume={33},
  number={5},
  pages={2223--2273},
  year={2020},
  publisher={Oxford University Press}
}

@article{markowitz1952portfolio,
  title={Portfolio Selection},
  author={Markowitz, HM},
  journal={The Journal of Finance},
  volume={1},
  pages={71--91},
  year={1952}
}

@article{rockafellar2000optimization,
  title={Optimization of conditional value-at-risk},
  author={Rockafellar, R Tyrrell and Uryasev, Stanislav and others},
  journal={Journal of risk},
  volume={2},
  pages={21--42},
  year={2000},
  publisher={Citeseer}
}

@article{davis1990portfolio,
  title={Portfolio selection with transaction costs},
  author={Davis, Mark HA and Norman, Andrew R},
  journal={Mathematics of operations research},
  volume={15},
  number={4},
  pages={676--713},
  year={1990},
  publisher={INFORMS}
}

@article{minh2018deep,
  title={Deep learning approach for short-term stock trends prediction based on two-stream gated recurrent unit network},
  author={Minh, Dang Lien and Sadeghi-Niaraki, Abolghasem and Huy, Huynh Duc and Min, Kyungbok and Moon, Hyeonjoon},
  journal={Ieee Access},
  volume={6},
  pages={55392--55404},
  year={2018},
  publisher={IEEE}
}

@article{song2019study,
  title={A study on novel filtering and relationship between input-features and target-vectors in a deep learning model for stock price prediction},
  author={Song, Yoojeong and Lee, Jae Won and Lee, Jongwoo},
  journal={Applied Intelligence},
  volume={49},
  pages={897--911},
  year={2019},
  publisher={Springer}
}

@article{go2019prediction,
  title={Prediction of stock value using pattern matching algorithm based on deep learning},
  author={Go, Yoon Ho and Hong, Jin Keun},
  journal={International Journal of Recent Technology and Engineering},
  volume={8},
  number={2},
  pages={31--35},
  year={2019}
}

@article{aggarwal2017deep,
  title={Deep investment in financial markets using deep learning models},
  author={Aggarwal, Saurabh and Aggarwal, Somya},
  journal={International Journal of Computer Applications},
  volume={162},
  number={2},
  pages={40--43},
  year={2017},
  publisher={Foundation of Computer Science}
}

@article{fang2019research,
  title={Research on quantitative investment strategies based on deep learning},
  author={Fang, Yujie and Chen, Juan and Xue, Zhengxuan},
  journal={Algorithms},
  volume={12},
  number={2},
  pages={35},
  year={2019},
  publisher={MDPI}
}

@book{markowitz2000mean,
  title={Mean-variance analysis in portfolio choice and capital markets},
  author={Markowitz, Harry M and Todd, G Peter},
  volume={66},
  number={4},
  year={2000},
  publisher={John Wiley \& Sons}
}

@book{jorion1997value,
  title={Value at Risk: The New Benchmark for Controlling Market Risk},
  author={Jorion, P.},
  isbn={9780786308484},
  lccn={96021381},
  url={https://books.google.co.in/books?id=u8efQgAACAAJ},
  year={1997},
  publisher={McGraw-Hill}
}

@article{konno1993mean,
  title={A mean-absolute deviation-skewness portfolio optimization model},
  author={Konno, Hiroshi and Shirakawa, Hiroshi and Yamazaki, Hiroaki},
  journal={Annals of Operations Research},
  volume={45},
  number={1},
  pages={205--220},
  year={1993},
  publisher={Springer}
}

@article{young1998minimax,
  title={A minimax portfolio selection rule with linear programming solution},
  author={Young, Martin R},
  journal={Management science},
  volume={44},
  number={5},
  pages={673--683},
  year={1998},
  publisher={INFORMS}
}

@article{nawrocki1992characteristics,
  title={The characteristics of portfolios selected by n-degree lower partial moment},
  author={Nawrocki, David N},
  journal={International Review of Financial Analysis},
  volume={1},
  number={3},
  pages={195--209},
  year={1992},
  publisher={Elsevier}
}

@article{brogan2008non,
  title={Non-separation in the mean--lower-partial-moment portfolio optimization problem},
  author={Brogan, Anita J and Stidham Jr, Shaler},
  journal={European Journal of Operational Research},
  volume={184},
  number={2},
  pages={701--710},
  year={2008},
  publisher={Elsevier}
}

@article{chaouki2020deep,
  title={Deep deterministic portfolio optimization},
  author={Chaouki, Ayman and Hardiman, Stephen and Schmidt, Christian and S{\'e}ri{\'e}, Emmanuel and De Lataillade, Joachim},
  journal={The Journal of Finance and Data Science},
  volume={6},
  pages={16--30},
  year={2020},
  publisher={Elsevier}
}

@article{moody2001learning,
  title={Learning to trade via direct reinforcement},
  author={Moody, John and Saffell, Matthew},
  journal={IEEE transactions on neural Networks},
  volume={12},
  number={4},
  pages={875--889},
  year={2001},
  publisher={IEEE}
}

@article{aboussalah2020continuous,
  title={Continuous control with stacked deep dynamic recurrent reinforcement learning for portfolio optimization},
  author={Aboussalah, Amine Mohamed and Lee, Chi Guhn},
  journal={Expert Systems with Applications},
  volume={140},
  pages={112891},
  year={2020},
  publisher={Elsevier}
}

@article{samuelson1970fundamental,
  title={The fundamental approximation theorem of portfolio analysis in terms of means, variances and higher moments},
  author={Samuelson, Paul A},
  journal={The Review of Economic Studies},
  volume={37},
  number={4},
  pages={537--542},
  year={1970},
  publisher={Wiley-Blackwell}
}

@article{wang2018leveraging,
  title={Leveraging deep learning with LDA-based text analytics to detect automobile insurance fraud},
  author={Wang, Yibo and Xu, Wei},
  journal={Decision Support Systems},
  volume={105},
  pages={87--95},
  year={2018},
  publisher={Elsevier}
}

@inproceedings{han2018nextgen,
  title={Distributed deep learning based language technologies to augment anti money laundering investigation},
  author={Han, Jingguang and Barman, Utsab and Hayes, Jer and Du, Jinhua and Burgin, Edward and Wan, Dadong},
   pages={37--42},
  year={2018},
  organization={Association for Computational Linguistics}
}

@article{smalter2017macroeconomic,
  title={Macroeconomic indicator forecasting with deep neural networks},
  author={Smalter Hall, Aaron and Cook, Thomas R},
  journal={Federal Reserve Bank of Kansas City Working Paper},
  number={17-11},
  year={2017}
}

@article{addo2018credit,
  title={Credit risk analysis using machine and deep learning models},
  author={Addo, Peter Martey and Guegan, Dominique and Hassani, Bertrand},
  journal={Risks},
  volume={6},
  number={2},
  pages={38},
  year={2018},
  publisher={MDPI}
}

@article{loureiro2018exploring,
  title={Exploring the use of deep neural networks for sales forecasting in fashion retail},
  author={Loureiro, Ana LD and Migu{\'e}is, Vera L and Da Silva, Lucas FM},
  journal={Decision Support Systems},
  volume={114},
  pages={81--93},
  year={2018},
  publisher={Elsevier}
}

@article{fombellida2020tackling,
  title={Tackling business intelligence with bioinspired deep learning},
  author={Fombellida, Juan and Mart{\'\i}n Rubio, Irene and Torres Alegre, Santiago and Andina, Diego},
  journal={Neural Computing and Applications},
  volume={32},
  pages={13195--13202},
  year={2020},
  publisher={Springer}
}

@INPROCEEDINGS{liang2013large,
  author={Liang, J. J. and Qu, B. Y.},
  booktitle={2013 IEEE Symposium on Swarm Intelligence (SIS)}, 
  title={Large-scale portfolio optimization using multiobjective dynamic mutli-swarm particle swarm optimizer}, 
  year={2013},
  volume={},
  number={},
  pages={1-6}}

@article{alvarez2023optimal,
  title={Optimal Brokerage Contracts in Almgren--Chriss Model with Multiple Clients},
  author={Alvarez, Guillermo Alonso and Nadtochiy, Sergey and Webster, Kevin},
  journal={SIAM Journal on Financial Mathematics},
  volume={14},
  number={3},
  pages={855--878},
  year={2023},
  publisher={SIAM}
}

@article{marcozzi2008stochastic,
  title={Stochastic optimal control of ultradiffusion processes with application to dynamic portfolio management},
  author={Marcozzi, Michael D},
  journal={Journal of computational and applied mathematics},
  volume={222},
  number={1},
  pages={112--127},
  year={2008},
  publisher={Elsevier}
}

@article{ekeland2007optimal,
  title={Optimal bond portfolios},
  author={Ekeland, Ivar and Taflin, Erik},
  journal={Paris-Princeton Lectures on Mathematical Finance },
  volume={1919},
  pages={51--102},
  year={2007},
  publisher={Springer}
}

@article{nadtochiy2019optimal,
  title={Optimal contract for a fund manager with capital injections and endogenous trading constraints},
  author={Nadtochiy, Sergey and Zariphopoulou, Thaleia},
  journal={SIAM Journal on Financial Mathematics},
  volume={10},
  number={3},
  pages={698--722},
  year={2019},
  publisher={SIAM}
}

@article{hambly2017stochastic,
  title={Stochastic evolution equations for large portfolios of stochastic volatility models},
  author={Hambly, Ben and Kolliopoulos, Nikolaos},
  journal={SIAM Journal on Financial Mathematics},
  volume={8},
  number={1},
  pages={962--1014},
  year={2017},
  publisher={SIAM}
}

@article{duedahl2016implementation,
  title={Implementation of Stochastic Yield Curve Duration and Portfolio Immunization Strategies},
  author={Duedahl, Sindre},
  journal={Journal of Mathematical Finance},
  volume={6},
  number={3},
  pages={401--415},
  year={2016},
  publisher={Scientific Research Publishing}
}

@article{pesenti2023portfolio,
  title={Portfolio optimization within a Wasserstein ball},
  author={Pesenti, Silvana M and Jaimungal, Sebastian},
  journal={SIAM Journal on Financial Mathematics},
  volume={14},
  number={4},
  pages={1175--1214},
  year={2023},
  publisher={SIAM}
}

@article{udeani2021application,
  title={Application of maximal monotone operator method for solving Hamilton--Jacobi--Bellman equation arising from optimal portfolio selection problem},
  author={Udeani, Cyril Izuchukwu and {\v{S}}ev{\v{c}}ovi{\v{c}}, Daniel},
  journal={Japan Journal of Industrial and Applied Mathematics},
  volume={38},
  number={3},
  pages={693--713},
  year={2021},
  publisher={Springer}
}

@article{bensoussan2022inter,
  title={Inter-temporal mutual-fund management},
  author={Bensoussan, Alain and Cheung, Ka Chun and Li, Yiqun and Yam, Sheung Chi Phillip},
  journal={Mathematical Finance},
  volume={32},
  number={3},
  pages={825--877},
  year={2022},
  publisher={Wiley Online Library}
}

@article{bensoussan2010real,
  title={Real options games in complete and incomplete markets with several decision makers},
  author={Bensoussan, Alain and Diltz, J David and Hoe, SingRu},
  journal={SIAM Journal on Financial Mathematics},
  volume={1},
  number={1},
  pages={666--728},
  year={2010},
  publisher={SIAM}
}

@article{miller2017optimal,
  title={Optimal control of conditional value-at-risk in continuous time},
  author={Miller, Christopher W and Yang, Insoon},
  journal={SIAM Journal on Control and Optimization},
  volume={55},
  number={2},
  pages={856--884},
  year={2017},
  publisher={SIAM}
}

@article{gatta2023meshless,
  title={Meshless methods for American option pricing through physics-informed neural networks},
  author={Gatta, Federico and Di Cola, Vincenzo Schiano and Giampaolo, Fabio and Piccialli, Francesco and Cuomo, Salvatore},
  journal={Engineering Analysis with Boundary Elements},
  volume={151},
  pages={68--82},
  year={2023},
  publisher={Elsevier}
}

@article{kim2024physics,
  title={Physics-informed convolutional transformer for predicting volatility surface},
  author={Kim, Soohan and Yun, Seok Bae and Bae, Hyeong Ohk and Lee, Muhyun and Hong, Youngjoon},
  journal={Quantitative Finance},
  volume={24},
  number={2},
  pages={203--220},
  year={2024},
  publisher={Taylor \& Francis}
}

@article{wang2023deep,
  title={A deep learning based numerical PDE method for option pricing},
  author={Wang, Xiang and Li, Jessica and Li, Jichun},
  journal={Computational economics},
  volume={62},
  number={1},
  pages={149--164},
  year={2023},
  publisher={Springer}
}

@article{bai2022application,
  title={The application of improved physics-informed neural network  method in finance},
  author={Bai, Yuexing and Chaolu, Temuer and Bilige, Sudao},
  journal={Nonlinear Dynamics},
  volume={107},
  number={4},
  pages={3655--3667},
  year={2022},
  publisher={Springer}
}

@article{carta2021multi,
  title={A multi-layer and multi-ensemble stock trader using deep learning and deep reinforcement learning},
  author={Carta, Salvatore and Corriga, Andrea and Ferreira, Anselmo and Podda, Alessandro Sebastian and Recupero, Diego Reforgiato},
  journal={Applied Intelligence},
  volume={51},
  pages={889--905},
  year={2021},
  publisher={Springer}
}

@article{gueant2019deep,
  title={Deep reinforcement learning for market making in corporate bonds: beating the curse of dimensionality},
  author={Gu{\'e}ant, Olivier and Manziuk, Iuliia},
  journal={Applied Mathematical Finance},
  volume={26},
  number={5},
  pages={387--452},
  year={2019},
  publisher={Taylor \& Francis}
}

@article{deng2016deep,
  title={Deep direct reinforcement learning for financial signal representation and trading},
  author={Deng, Yue and Bao, Feng and Kong, Youyong and Ren, Zhiquan and Dai, Qionghai},
  journal={IEEE transactions on neural networks and learning systems},
  volume={28},
  number={3},
  pages={653--664},
  year={2016},
  publisher={IEEE}
}

@article{hu2021survey,
  title={A survey of forex and stock price prediction using deep learning},
  author={Hu, Zexin and Zhao, Yiqi and Khushi, Matloob},
  journal={Applied System Innovation},
  volume={4},
  number={1},
  pages={9},
  year={2021},
  publisher={Mdpi}
}

@article{wang2020portfolio1,
  title={Portfolio formation with preselection using deep learning from long-term financial data},
  author={Wang, Wuyu and Li, Weizi and Zhang, Ning and Liu, Kecheng},
  journal={Expert Systems with Applications},
  volume={143},
  pages={113042},
  year={2020},
  publisher={Elsevier}
}

@article{tsang2020deep,
  title={Deep-learning solution to portfolio selection with serially dependent returns},
  author={Tsang, Ka Ho and Wong, Hoi Ying},
  journal={SIAM Journal on Financial Mathematics},
  volume={11},
  number={2},
  pages={593--619},
  year={2020},
  publisher={SIAM}
}

@article{soleymani2021deep,
  title={Deep graph convolutional reinforcement learning for financial portfolio management--DeepPocket},
  author={Soleymani, Farzan and Paquet, Eric},
  journal={Expert Systems with Applications},
  volume={182},
  pages={115127},
  year={2021},
  publisher={Elsevier}
}

@article{raissi2019physics,
  title={Physics-informed neural networks: A deep learning framework for solving forward and inverse problems involving nonlinear partial differential equations},
  author={Raissi, Maziar and Perdikaris, Paris and Karniadakis, George E},
  journal={Journal of Computational physics},
  volume={378},
  pages={686--707},
  year={2019},
  publisher={Elsevier}
}
\end{document}